\documentclass[10pt,twocolumn,twoside]{IEEEtran}

\usepackage{bm}   
\usepackage{amsmath}
\usepackage{amssymb}
\usepackage{indentfirst}  
\usepackage{graphicx}  
\usepackage{subfigure}
\usepackage{multirow} 
\usepackage{booktabs} 
\usepackage{cite}  
\usepackage{threeparttable}   
\usepackage{algorithm}
\usepackage{algpseudocode}
\usepackage{arydshln}
\usepackage{color}
\usepackage[american]{babel}
\usepackage{microtype}
\usepackage[numbers,sort&compress]{natbib}  

\usepackage{cases}
\usepackage{enumerate}

\newtheorem{theorem}{Theorem}
\newtheorem{corollary}{Corollary}
\newtheorem{lemma}{Lemma}
\newtheorem{definition}{Definition}
\newtheorem{remark}{Remark}
\newtheorem{assumption}{Assumption}

\def\degree{${}^{\circ}$}    
\def\mV{\mathcal{V}}
\def\mE{\mathcal{E}}

\def\mG{\mathcal{G}}

\def\mC{\mathcal{C}}

\def\mN{\mathcal{N}}

\def\mR{\mathcal{R}}

\def\mbR{\mathbb{R}}
\def\mbC{\mathbb{C}}
\def\mP{\mathcal{P}}
\def\mR{\mathcal{R}}

\def\mH{\mathcal{H}}

\def\Lg{\bm{L}_{\mG}}
\def\Lgs{\bm{L}_{\mG_*}}
\def\Lgz{\bm{L}_{\mG_0}}
\def\Lgr{\bm{L}_{\mG_r}}

\def\reff{r_{\textup{eff}}}
\def\geff{g_{\textup{eff}}}
\def\upj{\textup{j}}
\newcommand{\diag}[1]{\textup{diag}\{#1\}}

\begin{document}
%
\title{On Extension of Effective Resistance with Application to Graph Laplacian Definiteness \\and Power Network Stability}

\author{Yue~Song,~\IEEEmembership{Member,~IEEE,}
        David~J.~Hill,~\IEEEmembership{Life Fellow,~IEEE,}
        and~Tao~Liu,~\IEEEmembership{Member,~IEEE}


\thanks{Y. Song and T. Liu are with the Department of Electrical and Electronic Engineering, The University of Hong Kong, Hong Kong (e-mail: yuesong@eee.hku.hk; taoliu@eee.hku.hk).}
\thanks{D. J. Hill is with the Department of Electrical and Electronic Engineering, The University of Hong Kong, Hong Kong, and also with the School of Electrical and Information Engineering, The University of Sydney, Sydney, NSW 2006, Australia (e-mail: dhill@eee.hku.hk; david.hill@sydney.edu.au).}
}

\markboth{IEEE TRANSACTIONS ON CIRCUITS AND SYSTEMS--I: REGULAR PAPERS (TO APPEAR)} 
{Song \MakeLowercase{\textit{et al.}}: On Extension of Effective Resistance}

\maketitle

\begin{abstract}
  This paper extends the definitions of effective resistance and effective conductance to characterize the overall relation (positive coupling or antagonism) between any two disjoint sets of nodes in a signed graph. It generalizes the traditional definitions that only apply to a pair of nodes. The monotonicity and convexity properties are preserved by the extended definitions.
  The extended definitions provide new insights into graph Laplacian definiteness and power network stability.
  It is proved that the Laplacian matrix of a signed graph is positive semi-definite with only one zero eigenvalue if and only if the effective conductances between some specific pairs of node sets are positive. Also the number of Laplacian negative eigenvalues is upper bounded by the number of negative weighted edges.
  In addition, new conditions for the small-disturbance angle stability, hyperbolicity and type of power system equilibria are established, which intuitively interpret angle instability as the electrical antagonism between certain two sets of nodes in the defined active power flow graph.
  Moreover, a novel optimal power flow (OPF) model with effective conductance constraints is formulated, which significantly enhances power system transient stability.
  By the properties of extended effective conductance, the proposed OPF model admits a convex relaxation representation that achieves global optimality.
\end{abstract}


%
\IEEEpeerreviewmaketitle


\section{Introduction}\label{secintro}
The effective resistance and effective conductance\footnote{The effective resistance and effective conductance will be used interchangeably in this paper as they are reciprocal to each other.} are useful graph theory concepts originating from electrical networks \cite{klein1993resistance}, the history of which can date back to 1940s \cite{foster1949average}.
In the context of resistive networks, the effective resistance describes the port resistance between a pair of nodes that captures the global properties of network, i.e., both direct connections and indirect connections between the two nodes.
In addition, the mathematical foundations of effective resistance, such as distance metric, monotonicity and convexity, have been well established \cite{bapat2010graphs, ghosh2008minimizing}.

\indent
Due to the clear physical meaning and useful properties, the effective resistance has widespread applications in theoretical and practical problems.
For instance, the effective resistance is closely connected to other graph theory concepts such as Kirchhoff Index, random walks and Foster's theorem \cite{thulasiraman2019network}, which have applications in, e.g., designing online algorithms \cite{coppersmith1993random} and describing molecular structure in chemistry \cite{bonchev1994molecular}.
The effective resistance has also been used in many fields of electric power networks, such as cascading failures \cite{guo2017monotonicity, soltan2017analysis}, network partitioning \cite{cotilla2013multi} and power network stability \cite{dorfler2010spectral, dorfler2018electrical, song2017networkbased}.
In addition, since the effective resistance is defined by the graph Laplacian matrix \cite{klein1993resistance}, it is useful in quantifying the performance of network control problems where the graph Laplacian plays an important role.
Some bounds for the Laplacian spectrum are constructed by effective resistance in \cite{barooah2006graph}, which describe the convergence speed of distributed control.
The $\mH_2$ norm and $\mH_{\infty}$ norm of linear oscillator networks are given in terms of effective resistance in \cite{grunberg2018performance} and \cite{johansson2018optimization}, respectively.
In \cite{young2016new1, young2016new2}, the definition of effective resistance is extended from undirected graphs to directed graphs, which provides a new tool for control problems over directed graphs.

\indent
Recently, the effective resistance concept has been introduced to signed graphs that have both positive and negative weighted edges.
The negative weighted edges can model the attacks on network topology or antagonistic effects between some nodes.
The Laplacian matrix of a signed graph (or simply called the signed Laplacian in some literature \cite{chen2016characterizing}) being positive semi-definite (PSD) with only one zero eigenvalue is a key requirement for reaching network consensus or synchronization.
It is revealed in \cite{zelazo2014definiteness, zelazo2017robustness, chen2016definiteness, ahmadizadeh2017eigenvalues} that the sign of effective resistance can be used to check the definiteness of Laplacian matrix, which applies to a special class of signed graphs where the negative weighted edges are not contained in the same cycle.
By using the concept of mutual effective resistance between two node pairs, a linear matrix inequality condition for Laplacian definiteness is proposed in \cite{chen2016characterizing, chen2017spectral}, which applies to generic signed graphs.
Nevertheless, there is still a lack of a concise scalar condition for the Laplacian definiteness of generic signed graphs.

\indent
So far the effective resistance has been defined in the ``node-to-node'' manner that describes the relation between two nodes.
However, the collective behaviors observed in many practical problems over networks are concerned with the ``nodes-to-nodes'' relation.
For instance, the instability event in power networks occurs in such a form that a group of generators go unstable with respect to the remaining ones \cite{song2017cutset}.
By the actual requirements observed from practical problems and the aforementioned bottlenecks in characterizing Laplacian definiteness, the node-to-node description of effective resistance could be a limitation.
On the other hand, a nodes-to-nodes version of effective resistance, which describes the relation between two sets of nodes, is a natural direction for extension that could offer new thoughts into these issues.

\indent
Following this idea, we extend the definitions of effective resistance and effective conductance to the ``nodes-to-nodes'' version, and apply them to the study of graph Laplacian definiteness and power network stability.
The contributions of this paper are threefold:

\indent
1) The extended definitions generalize the traditional ones with intuitive circuit interpretations from the perspectives of current flow, energy dissipation and potential difference.
The extended effective resistance/conductance can be used to characterize the overall coupling between any two disjoint sets of nodes.
When applied to signed graphs, the positive or negative sign of the extended effective resistance/conductance indicates the positive coupling or antagonism between the two sets of nodes.
Also the monotonicity and convexity properties are preserved by the extended definitions.

\indent
2) The extended definitions shed new light on graph Laplacian definiteness.
We establish some necessary and sufficient conditions for the Laplacian matrix of a generic signed graph being PSD with only one zero eigenvalue, which links the Laplacian definiteness to the positive effective conductances between specific pairs of node sets.
By the properties of extended effective conductance, we also prove that the number of Laplacian negative eigenvalues is upper bounded by the number of negative weighted edges.

\indent
3) The extended definitions shed new light on power network stability.
We establish some necessary and sufficient conditions for the small-disturbance angle stability, hyperbolicity and type of an equilibrium point in terms of the extended effective conductance.
These conditions generalize our previous results in \cite{song2017networkbased}, by which the instability mechanism can be intuitively interpreted as the electrical antagonism between certain two sets of nodes.
A new optimal power flow (OPF) model is also developed, which enhances transient stability by enforcing the effective conductance constraint with respect to coherent generators.
The OPF model admits a convex relaxation formulation due to the useful properties of extended effective conductance, which promotes global optimality.

\indent
The rest of the paper is organized as follows.
Some background knowledge and a brief review of the traditional effective resistance and effective conductance are given in Section \ref{secformu}.
In Section \ref{secextend}, the extended definitions of effective resistance and effective conductance are proposed and some basic properties are studied. In Section \ref{seclapPSD} and Section \ref{secstability}, the extended definitions are applied to graph Laplacians and power networks, respectively; some new results on graph Laplacian definiteness and power network stability are established.
Section \ref{secconclu} concludes the paper and makes future prospects.
Those proofs with tedious derivation are given in Appendix for better readability.

\section{Background and problem setup}\label{secformu}
\subsection{Preliminaries and notations}
For simplicity, we use $\bm{x}=[x_i]\in \mbR^{p}$ to represent a vector, and $\bm{x}=\diag{x_i}\in \mbR^{p\times p}$ to represent a diagonal matrix.
The triple $\{i_+(\bm{A}), i_-(\bm{A}), i_0(\bm{A})\}$ denotes the matrix inertia of $\bm{A}\in \mbR^{p\times p}$, i.e., the number of eigenvalues with positive, negative and zero parts, respectively.
$\bm{A}\succ 0$ (or $\bm{A}\succeq 0$) means $\bm{A}$ is positive definite (or positive semi-definite).
The notation $\bm{I}_p\in \mbR^{p\times p}$ denotes an identity matrix, $\bm{1}_{p}\in \mbR^{p}$ denotes a vector with all entries being one, $\bm{0}_{p}\in \mbR^{p}$ denotes the zero vector, $\bm{A}^{\dag}$ denotes the Moore-Penrose inverse of $\bm{A}$, and $|\cdot|$ denotes the cardinality of a set.
The italic $j$ denotes a numbering index, and the upright j denotes the square root of -1.
The notation $\bm{e}_i$ denotes a vector with its $i$-th entry being one and the other entries being zero, $\bm{e}_{\mV_a}$ denotes a vector with those entries indexed by the set $\mV_a$ being one and the other entries being zero. Note that the dimensions of $\bm{e}_i,\bm{e}_{\mV_a}$ are not fixed in the paper, they will be consistent with other matrices in the respective equations.
Let $\bm{H} =
      \begin{bmatrix}
        \bm{A} & \bm{B} \\
        \bm{C} & \bm{D} \\
      \end{bmatrix}$
be a matrix with two-by-two block partition, where $\bm{A}\in \mbR^{p\times p}$, $\bm{B}\in \mbR^{p\times q}$, $\bm{C}\in \mbR^{q\times p}$ and $\bm{D}\in \mbR^{q\times q}$.
Assume $\bm{D}$ is nonsingular, the Schur complement of the block $\bm{D}$ of the matrix $\bm{H}$ is denoted as $\bm{H}/ \bm{D}=\bm{A} - \bm{B}\bm{D}^{-1}\bm{C}$.

\indent
Denote $\mG(\mV,\mE,\bm{W})$ as a weighted undirected graph, where $\mV$ is the set of nodes with $n=|\mV|$, $\mE\subseteq\mV\times\mV$ is the set of edges with $l=|\mE|$, and $\bm{W}\in \mbR^{l\times l}$ is the diagonal matrix indicating edge weights.
The edge $k$ connecting node $i$ and node $j$ is denoted by either $\epsilon_k$ or an unordered pair $(i,j)$ (i.e., $(i,j)$ and $(j,i)$ are equivalent notations).
The weight of edge $\epsilon_k=(i,j)\in \mE$ is given by the $k$-th main diagonal of $\bm{W}$, denoted $w_{ij}$. We consider signed graphs where each edge can have a positive or negative weight.
For defining the incidence matrix, suppose each edge is fictitiously assigned an arbitrary but fixed orientation. Then, the incidence matrix $\bm{E}\in \mbR^{n\times l}$ is defined such that $\forall \epsilon_k=(i,j)\in\mE$, $E_{ik}=1$ if $\epsilon_k$ originates at node $i$, $E_{jk}=-1$ if $\epsilon_k$ terminates at node $j$ and $E_{mk}=0, m\neq i,j$.
The Laplacian matrix is defined as $\Lg=\bm{EWE}^{T}\in \mbR^{n\times n}$, and we use $\bm{L}_{xy}$ to denote the submatrix of $\Lg$ with the rows indexed by the set of nodes $\mV_x\subseteq\mV$ and columns indexed by the set of nodes $\mV_y\subseteq\mV$.
Note that $\Lg$ has at least one zero eigenvalue as $\Lg\bm{1}_n=\bm{0}_n$.
Without loss of generality, the graphs studied in this paper are assumed to be connected, while the obtained results also apply to each connected component of a disconnected graph.

\indent
In addition, we present two lemmas below that will be frequently used in the paper.

\begin{lemma}\label{lemLdag}\cite{gutman2004generalized}
   If $\Lg$ is PSD with only one zero eigenvalue, then $\Lg^{\dag}\Lg=\Lg\Lg^{\dag}=\bm{I}_n-\frac{1}{n}\bm{1}_n\bm{1}_n^T$.
\end{lemma}

\begin{lemma}\label{lemschur}\cite{zhang2006schur}
   Let $\bm{H} =
      \begin{bmatrix}
        \bm{A} & \bm{B} \\
        \bm{B}^{T} & \bm{D} \\
      \end{bmatrix}$
   be a real symmetric matrix with the square block $\bm{D}$ being nonsingular.
   Then $i_+(\bm{H})=i_+(\bm{H}/ \bm{D}) + i_+(\bm{D})$, $i_-(\bm{H})=i_-(\bm{H}/ \bm{D}) + i_-(\bm{D})$ and $i_0(\bm{H})=i_0(\bm{H}/ \bm{D})$.
\end{lemma}

\subsection{Effective resistance/conductance: a brief review}
Let us have a brief review of the traditional effective resistance and effective conductance.
The currently common expression for effective resistance is proposed in \cite{klein1993resistance}.
The reciprocal of effective resistance is referred to as effective conductance \cite{doyle1984random, dorfler2010spectral}.
The detailed definitions are given below.

\begin{definition}[Effective Resistance/Conductance \cite{klein1993resistance, doyle1984random}]
   For any pair of nodes $i,j\in \mV$, $i\neq j$, the effective resistance between $i$ and $j$ is defined as
   \begin{equation}\label{Rij}
   \begin{split}
      \mR_{ij}=(\bm{e}_i-\bm{e}_j)^T\Lg^{\dag}(\bm{e}_i-\bm{e}_j)
   \end{split}
   \end{equation}
   and the effective conductance between $i$ and $j$ is defined as
   \begin{equation}\label{Cij}
   \begin{split}
      \mC_{ij}=\frac{1}{(\bm{e}_i-\bm{e}_j)^T\Lg^{\dag}(\bm{e}_i-\bm{e}_j)}.
   \end{split}
   \end{equation}
   where $\bm{e}_i, \bm{e}_j\in\mbR^n$ in \eqref{Rij} and \eqref{Cij}.
\end{definition}

\begin{figure}[!h]
  \centering
  \includegraphics[width=3.5in]{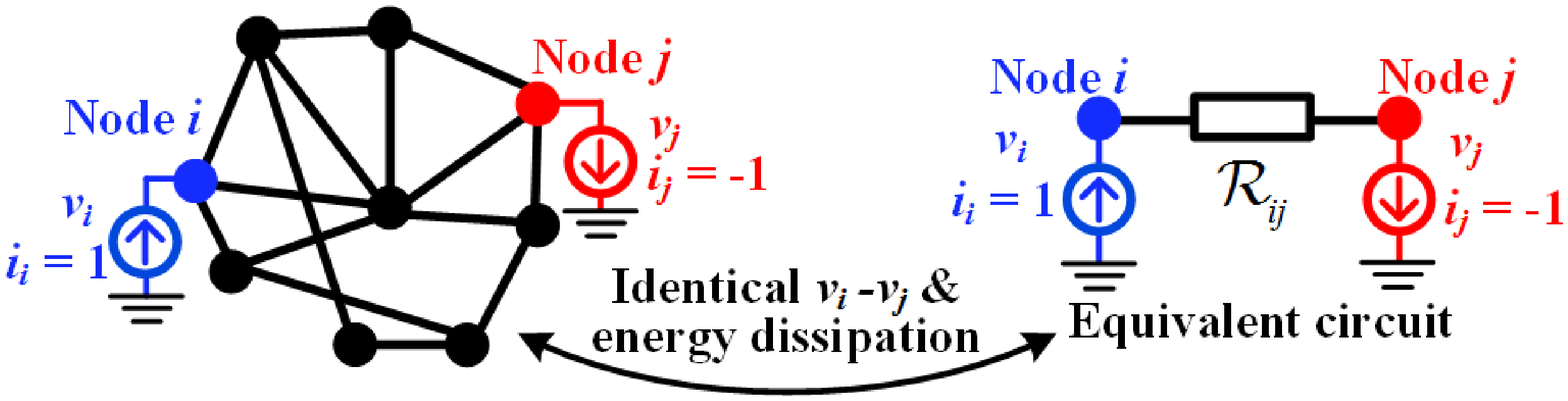}
  \caption{Circuit interpretation of $\mR_{ij}$.}
  \label{figgefftraditional}
\end{figure}

\indent
The concept of effective resistance originates from electrical networks with clear physical meanings.
Consider a resistive network interpreted by the graph $\mG(\mV,\mE,\bm{W})$, i.e., it has the same topology as $\mG$ and the conductance of each transmission line $(i,j)$ is equal to $w_{ij}$.
Suppose node $i$ connects a source with unit current injection, node $j$ connects a source with unit current ejection and the other nodes keep open circuit. Then, applying Kirchhoff's circuit law gives $\bm{i} = \bm{e}_i-\bm{e}_j = \Lg\bm{v}$, where $\bm{i}, \bm{v}\in \mbR^n$ denote the current vector and voltage vector, respectively.
Suppose all lines have positive conductance, then $\Lg$ is PSD with only one zero eigenvalue \cite{ghosh2008minimizing}, and we have $\Lg^{\dag}\bm{i}=\bm{v}-\frac{1}{n}\bm{1}_n\bm{1}_n^T\bm{v}$ by Lemma \ref{lemLdag}.
Thus, the total energy dissipation of the network and potential difference between node $i$ and node $j$ are given by
\begin{equation}\label{PtotVdiff}
\begin{split}
      P_{\textup{tot}} &= \bm{i}^T\bm{v} = \bm{i}^T\Lg^{\dag}\bm{i} = \mR_{ij} \\
      v_i - v_j &= (\bm{e}_i-\bm{e}_j)^T\bm{v} = (\bm{e}_i-\bm{e}_j)^T\Lg^{\dag}\bm{i}  = \mR_{ij}.
\end{split}
\end{equation}
It implies that the port between node $i$ and node $j$ can be equivalent as a series resistance $\mR_{ij}$ (see Fig. \ref{figgefftraditional}). With the same current injection and ejection at node $i$ and node $j$, the equivalent circuit has the same energy dissipation and potential difference as the original network.
Hence, the effective resistance and effective conductance characterize the overall link between a pair of nodes. A smaller positive $\mR_{ij}$ (or equivalently, larger positive $\mC_{ij}$) indicates a stronger coupling between the pair of nodes.

\indent
As aforementioned, the effective resistance gains popularity in many subjects.
Nevertheless, the traditional definition has some limitations. It will be seen later that it may not work properly for signed graphs.
In addition, many engineering problems are concerned with the relation between two sets of nodes rather than a pair of nodes, e.g., a group of generators losing synchronism with respect to the remaining ones in case of power system transient instability \cite{song2017cutset}, where the traditional definition of effective resistance does not apply.
This appeals for an extension that will be addressed in the following.

\section{Extending effective resistance}\label{secextend}
\subsection{Extended definitions and physical interpretations}\label{secdefreff}
In this section, the concepts of effective resistance and effective conductance are extended as a proper characterization of the overall coupling between any two disjoint sets of nodes with the presence of negative weighted edges.
Given two non-empty disjoint sets of nodes, say $\mV_a, \mV_b\subset \mV$, and let $\mV_c=\mV\backslash (\mV_a\cup\mV_b)$ be the set of remaining nodes. For simplicity, denote $a=|\mV_a|$, $b=|\mV_b|$ and $c=|\mV_c|$.
Then $\Lg$ can be rewritten as
\begin{equation}
\begin{split}
   \Lg=
   \begin{bmatrix}
     \bm{L}_{aa}   &   \bm{L}_{ab}   &   \bm{L}_{ac}     \\
     \bm{L}_{ab}^T   &   \bm{L}_{bb}   &   \bm{L}_{bc}    \\
      \bm{L}_{ac}^T  &   \bm{L}_{bc}^T   &   \bm{L}_{cc}    \\
   \end{bmatrix}.
\end{split}
\end{equation}
Let $\Lg/ \bm{L}_{cc}$ denote the Schur complement of the block $\bm{L}_{cc}$ of the Laplacian matrix, which is real symmetric. In addition, we set $\Lg/ \bm{L}_{cc}=\Lg$ if $\mV_c=\phi$.
With these notations, we make an assumption and propose the extended definitions of effective resistance and effective conductance below.

\begin{assumption}\label{assump1}
   The matrix $\bm{L}_{cc}$ is nonsingular for any non-empty disjoint sets $\mV_a, \mV_b\subset \mV$ and $\mV_c=\mV\backslash (\mV_a\cup\mV_b)\neq \phi$.
\end{assumption}

\begin{definition}[Extended Effective Resistance/Conductance]\label{defreff}
   Let $\mV_a, \mV_b\subset \mV$ be two non-empty disjoint sets of nodes in the graph $\mG(\mV, \mE, \bm{W})$, and $\mV_c=\mV\backslash (\mV_a\cup\mV_b)$ be the set of remaining nodes.
   The extended effective resistance and effective conductance between $\mV_a$ and $\mV_b$, say $\reff(\Lg, \mV_a, \mV_b)$ and $\geff(\Lg, \mV_a, \mV_b)$, are respectively defined as
   \begin{subequations}
   \begin{align}
      \reff(\Lg, \mV_a, \mV_b)&=(\bm{e}_{\mV_a}^T(\Lg/ \bm{L}_{cc})\bm{e}_{\mV_a})^{-1} \label{equreff}  \\
      \geff(\Lg, \mV_a, \mV_b)&=\bm{e}_{\mV_a}^T(\Lg/ \bm{L}_{cc})\bm{e}_{\mV_a} \label{equgeff}
   \end{align}
   \end{subequations}
   where $\bm{e}_{\mV_a}\in\mbR^{a+b}$.
\end{definition}

\begin{figure}[!h]
  \centering
  \includegraphics[width=3.3in]{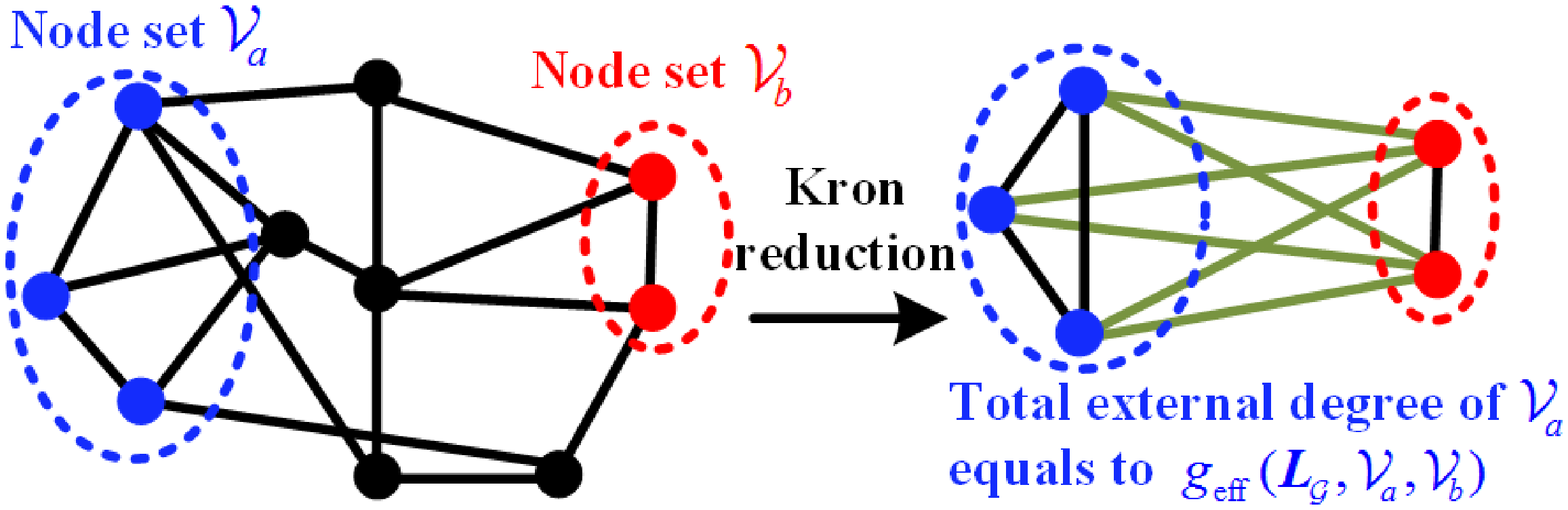}
  \caption{Direct meaning of $\geff(\Lg, \mV_a, \mV_b)$.}
  \label{figgeffdirect}
\end{figure}

\indent
The matrix $\Lg/ \bm{L}_{cc}$ is indeed the Laplacian matrix of the graph after eliminating $\mV_c$ via Kron reduction \cite{dorfler2013kron} (it will be shown in Lemma \ref{lemLcc} that $\Lg/ \bm{L}_{cc}$ has zero row-sums).
This observation leads to a straightforward meaning of $\geff(\Lg, \mV_a, \mV_b)$---the total external degree of $\mV_a$ in the reduced graph (e.g., the total weights of the green edges in Fig. \ref{figgeffdirect}).
In addition, it will be seen later that the extended definitions have other intuitive physical interpretations in the context of resistive networks, and they include the traditional definitions as a special case.
Moreover, we define $\reff(\Lg, \mV_a, \mV_b)=\infty$ in case that $\geff(\Lg, \mV_a, \mV_b)=0$, which corresponds to ``open-circuit'' status.

\indent
Note that Assumption \ref{assump1} is required by the well-definedness of $\reff(\Lg, \mV_a, \mV_b)$ and $\geff(\Lg, \mV_a, \mV_b)$.
All the following results will be based on this assumption, i.e., when we mention $\reff(\Lg, \mV_a, \mV_b)$ and $\geff(\Lg, \mV_a, \mV_b)$ in the lemmas, theorems and corollaries, it implicitly means the concerned $\bm{L}_{cc}$ is nonsingular.
The assumption can be justified by the fact that the singularity of $\bm{L}_{cc}$ is generally linked to nonphysical circuit modeling \cite{chua1980dynamic}.
For example, in the resistive network interpreted by $\mG(\mV, \mE, \bm{W})$, suppose each node $i\in \mV_a\cup\mV_b$ connects a shunt capacitor $C_i$, while each node $i\in \mV_c$ keeps open circuit. The network dynamics is then described by
\begin{equation*}
\begin{split}
   \begin{bmatrix}
      -\bm{C}_a\dot{\bm{v}}_a   \\
      -\bm{C}_b\dot{\bm{v}}_b   \\
      \bm{0}_c
   \end{bmatrix}=
   \begin{bmatrix}
     \bm{L}_{aa}   &   \bm{L}_{ab}   &   \bm{L}_{ac}     \\
     \bm{L}_{ab}^T   &   \bm{L}_{bb}   &   \bm{L}_{bc}    \\
      \bm{L}_{ac}^T  &   \bm{L}_{bc}^T   &   \bm{L}_{cc}    \\
   \end{bmatrix}
   \begin{bmatrix}
      \bm{v}_a   \\
      \bm{v}_b   \\
      \bm{v}_c
   \end{bmatrix}
\end{split}
\end{equation*}
where $\bm{v}_a, \bm{v}_b, \bm{v}_c$ denote the voltage vectors of $\mV_a, \mV_b, \mV_c$, respectively;
and $\bm{C}_a=\diag{C_i}\in \mbR^{a\times a}$, $\forall i\in \mV_a$ and $\bm{C}_b=\diag{C_i}\in \mbR^{b\times b}$, $\forall i\in \mV_b$. In case that $\bm{L}_{cc}$ is singular, the algebraic variable $\bm{v}_c$ becomes noncausal, i.e., it is not dependent on the state variable $\bm{v}_a, \bm{v}_b$.
In circuit theory, the network exhibits an impasse point and should be remodeled \cite{chua1980dynamic}. So it is reasonable to exclude such nonphysical situations in the analysis.
Moreover, when $\Lg$ is PSD with only one zero eigenvalue, it is trivial that all possible $\bm{L}_{cc}$, which refer to principal submatrices of $\Lg$, are nonsingular. Thus, for those results where $\Lg$ being PSD with only one zero eigenvalue is a precondition, Assumption \ref{assump1} is implicitly satisfied (e.g., see Theorem \ref{thmenergyVdiff}, Corollary \ref{corcompat}, Theorem \ref{thmmonotonenode}, Theorem \ref{thmmonotone}, Theorem \ref{thmconvexity} and necessity part of Corollary \ref{corPSD}).

\indent
We show some basic properties of $\geff(\Lg, \mV_a, \mV_b)$ before presenting more physical interpretations.

\begin{lemma}\label{lemLcc}
   Given a graph $\mG(\mV, \mE, \bm{W})$ and two non-empty disjoint sets $\mV_a, \mV_b\subset \mV$, and $\mV_c=\mV\backslash (\mV_a\cup\mV_b)$. The following equalities hold
   \begin{subequations}\label{Lcc}
   \begin{align}
      (\Lg/ \bm{L}_{cc})\bm{1}_{a+b}&=\bm{0}_{a+b}  \label{Lcc1}\\
      \geff(\Lg, \mV_a, \mV_b)&=\geff(\Lg, \mV_b, \mV_a)   \label{Lcc2}\\
      \geff(\Lg, \mV_a, \mV_b)&=\bm{1}_a^T(\bm{L}_{aa} - \bm{L}_{ac}\bm{L}_{cc}^{-1}\bm{L}_{ac}^T)\bm{1}_a \label{Lcc3aa} \\
      \geff(\Lg, \mV_a, \mV_b)&=-\bm{1}_a^T(\bm{L}_{ab} - \bm{L}_{ac}\bm{L}_{cc}^{-1}\bm{L}_{bc}^T)\bm{1}_b \label{Lcc3ab}\\
      \geff(\Lg, \mV_a, \mV_b)&=\bm{1}_b^T(\bm{L}_{bb} - \bm{L}_{bc}\bm{L}_{cc}^{-1}\bm{L}_{bc}^T)\bm{1}_b \label{Lcc3bb}
   \end{align}
   \end{subequations}
   where the matrices associated with $\mV_c$ in \eqref{Lcc3aa}-\eqref{Lcc3bb} are neglected in case that $\mV_c$ is empty.
\end{lemma}

\begin{IEEEproof}
   The proof is given in Appendix A.
\end{IEEEproof}

\begin{lemma}\label{lemquotient}
    Given a graph $\mG(\mV, \mE, \bm{W})$, a subset of nodes $\mV_r\subset \mV$ with $|\mV_r|\geq 2$, and $\mV_o=\mV\backslash\mV_r$.
    Suppose $\bm{L}_{oo}$ is nonsingular and let $\mG_r(\mV_r, \mE_r, \bm{W}_r)$ be the graph whose Laplacian matrix is $\Lg/ \bm{L}_{oo}$.
    Then $\geff(\Lg, \mV_a, \mV_b)=\geff(\Lgr, \mV_a, \mV_b)$ for any two non-empty disjoint sets $\mV_a, \mV_b\subset \mV_r$.
\end{lemma}

\begin{IEEEproof}
   The proof is given in Appendix A.
\end{IEEEproof}

\indent
Similar results to Lemma \ref{lemLcc} have been obtained in \cite{van2010characterization} for graphs with positive weighted edges only. Here we extend them to the case with negative weighted edges.
Recalling $\mR_{ij}=\mR_{ji}$ and $\mC_{ij}=\mC_{ji}$ held by the traditional definitions, Lemma \ref{lemLcc} implies that the extended definitions preserve this symmetry feature. So we do not need to specify the order of $\mV_a$ and $\mV_b$ when using the extended definitions.
In addition, Lemma \ref{lemquotient} implies that the values of $\reff(\Lg, \mV_a, \mV_b)$ and $\geff(\Lg, \mV_a, \mV_b)$ are preserved under Kron reduction.
These two lemmas will be used in the proofs of the following results.

\indent
We now come to the circuit interpretations of $\reff(\Lg, \mV_a, \mV_b)$ and $\geff(\Lg, \mV_a, \mV_b)$.

\begin{theorem}[Current Flow Interpretation]\label{thmIa}
    Given a resistive network interpreted by $\mG(\mV,\mE,\bm{W})$ and two non-empty disjoint sets $\mV_a, \mV_b\subset \mV$, $\geff(\Lg, \mV_a, \mV_b)$ is equal to the total current flow from $\mV_a$ to $\mV_b$ when each node $i\in\mV_a$ connects to a voltage source with unit potential, each node $i\in\mV_b$ is grounded with zero potential, and each node $i\in\mV_c=\mV\backslash (\mV_a\cup\mV_b)$ remains open circuit.
\end{theorem}

\begin{IEEEproof}
   It follows from the settings in the statement that
   \begin{equation}\label{currentflow}
   \begin{split}
   \begin{bmatrix}
      \bm{i}_a   \\
      \bm{i}_b   \\
      \bm{0}_c
   \end{bmatrix}=
   \begin{bmatrix}
     \bm{L}_{aa}   &   \bm{L}_{ab}   &   \bm{L}_{ac}     \\
     \bm{L}_{ab}^T   &   \bm{L}_{bb}   &   \bm{L}_{bc}    \\
      \bm{L}_{ac}^T  &   \bm{L}_{bc}^T   &   \bm{L}_{cc}    \\
   \end{bmatrix}
   \begin{bmatrix}
      \bm{1}_a   \\
      \bm{0}_b   \\
      \bm{v}_c
   \end{bmatrix}.
   \end{split}
   \end{equation}
   Reducing the third row in \eqref{currentflow} gives
   $\begin{bmatrix}
      \bm{i}_a   \\
      \bm{i}_b   \\
   \end{bmatrix}= (\Lg/\bm{L}_{cc})
   \begin{bmatrix}
      \bm{1}_a   \\
      \bm{0}_b   \\
   \end{bmatrix}$.
   So the total current flow from $\mV_a$ to $\mV_b$ is
   \begin{equation}\label{totalIa}
   \begin{split}
      \bm{1}_a^T\bm{i}_a =
   \begin{bmatrix}
      \bm{1}_a^T  &  \bm{0}_b^T   \\
   \end{bmatrix}
      (\Lg/\bm{L}_{cc})
      \begin{bmatrix}
      \bm{1}_a   \\
      \bm{0}_b  \\
   \end{bmatrix} = \geff(\Lg, \mV_a, \mV_b)
   \end{split}
   \end{equation}
   which completes the proof.
\end{IEEEproof}

\begin{theorem}[Energy Dissipation \& Potential Difference Interpretations]\label{thmenergyVdiff}
   Given a graph $\mG(\mV, \mE, \bm{W})$ and any two non-empty disjoint sets $\mV_a, \mV_b\subset \mV$.
   If $\Lg$ is PSD with only one zero eigenvalue, then we have $\reff(\Lg, \mV_a, \mV_b)>0$ and
   \begin{subequations}\label{energydis}
   \begin{align}
      \reff(\Lg, \mV_a, \mV_b)=\min_{\bm{i}}~&\bm{i}^T\Lg^{\dag}\bm{i} \label{energydis0} \\
      s.t.~&\bm{e}_{\mV_a}^T\bm{i} = 1  \label{energydis1}\\
        &\bm{e}_{\mV_b}^T\bm{i} = -1  \label{energydis2}\\
        &\bm{e}_i^T\bm{i} = 0,~\forall i\in \mV_c  \label{energydis3}
   \end{align}
   \end{subequations}
   where $\bm{e}_{\mV_a},\bm{e}_{\mV_b},\bm{e}_i (i\in\mV_c)\in\mbR^{n}$; and
   \begin{equation}\label{Vdiff}
   \begin{split}
      \reff(\Lg, \mV_a, \mV_b)=(\bm{e}_1-\bm{e}_2)^T(\bm{\mP}_{ab}^T\Lg\bm{\mP}_{ab})^{\dag}(\bm{e}_1-\bm{e}_2)
   \end{split}
   \end{equation}
   where $\bm{e}_1,\bm{e}_2\in\mbR^{c+2}$ and $\bm{\mP}_{ab} =\begin{bmatrix}
                   \bm{e}_{\mV_a} & \bm{e}_{\mV_b} & \bm{X}_c
                 \end{bmatrix}\in \mbR^{n\times (c+2)}$
   with $\bm{X}_c\in\mbR^{n\times c}$ collecting the column vectors $\bm{e}_i\in\mbR^n$, $\forall i\in\mV_c$.
\end{theorem}

\begin{IEEEproof}
    The proof is given in Appendix A.
\end{IEEEproof}

\begin{corollary}[Downward Compatibility]\label{corcompat}
   If $\Lg$ is PSD with only one zero eigenvalue, then $\reff(\Lg, \mV_a, \mV_b)=\mR_{ij}$ and $\geff(\Lg, \mV_a, \mV_b)=\mC_{ij}$ if $\mV_a=\{i\}$, $\mV_b=\{j\}$, $\forall i,j\in \mV$, $i\neq j$.
\end{corollary}

\begin{IEEEproof}
    If $\mV_a=\{i\}$, $\mV_b=\{j\}$, the vectors $\bm{e}_{\mV_a}, \bm{e}_{\mV_b}$ in \eqref{energydis} are reduced to $\bm{e}_{\mV_a}=\bm{e}_i$, $\bm{e}_{\mV_b}=\bm{e}_j$. Then, the optimum of problem \eqref{energydis} is $\bm{i}=\bm{e}_i-\bm{e}_j$, and by Theorem \ref{thmenergyVdiff} we have $\reff(\Lg, \{i\}, \{j\})=(\bm{e}_i-\bm{e}_j)^T\Lg^{\dag}(\bm{e}_i-\bm{e}_j)=\mR_{ij}$ and $\geff(\Lg, \{i\}, \{j\})=\mC_{ij}$.
\end{IEEEproof}

\indent
Theorem \ref{thmIa} provides a current flow description for $\geff(\Lg, \mV_a, \mV_b)$.
Theorem \ref{thmenergyVdiff} interprets $\reff(\Lg, \mV_a, \mV_b)$ from energy dissipation and potential difference viewpoints, which will be detailed in the following remarks.
Corollary \ref{corcompat} implies that the extended definitions include the traditional ones as a special case when $\mV_a=\{i\}$, $\mV_b=\{j\}$.
Note that Theorem \ref{thmIa} requires no precondition on $\Lg$, while the precondition of Theorem \ref{thmenergyVdiff} and Corollary \ref{corcompat} ($\Lg$ being PSD with only one zero eigenvalue) is held by normal circuits where all line conductances are positive.
So the physical interpretations indicated by these theorems have general applicability.

\begin{figure}[!h]
  \centering
  \includegraphics[width=3.3in]{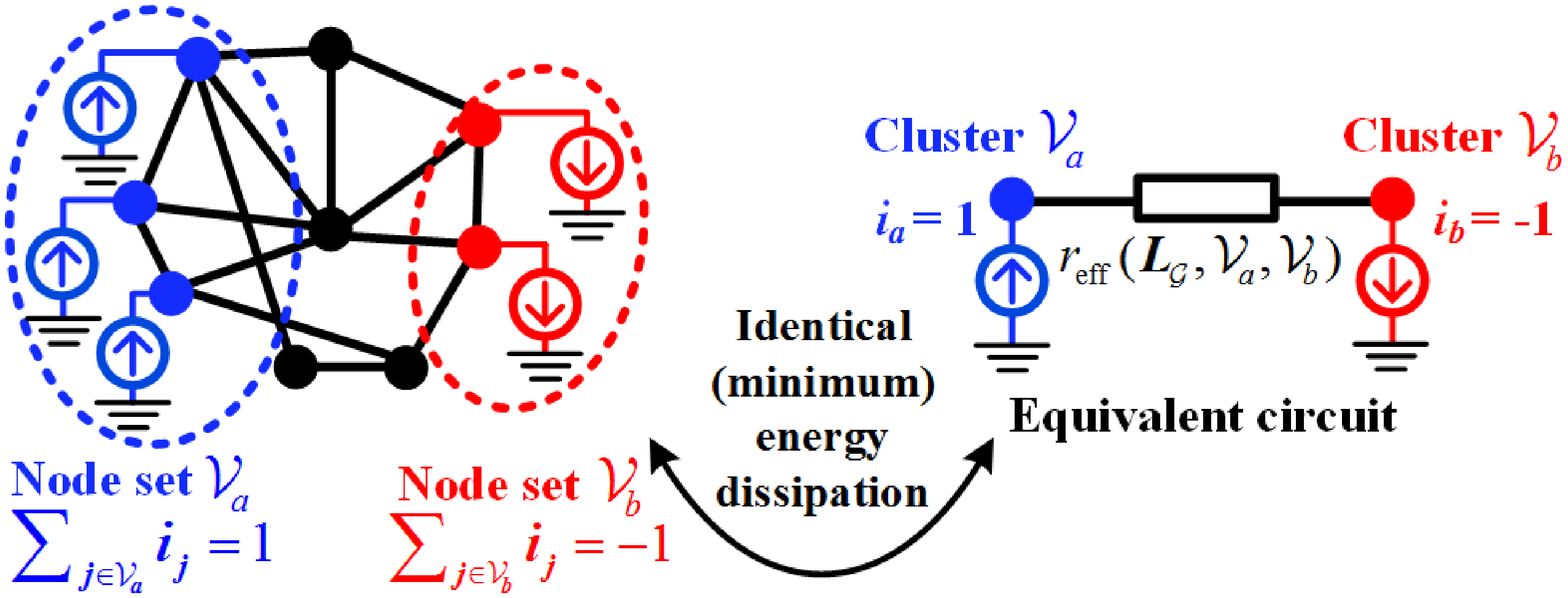}
  \caption{Interpretation of $\reff(\Lg, \mV_a, \mV_b)$ by energy dissipation.}
  \label{figgeffenergy}
\end{figure}

\begin{remark}
The optimization problem in \eqref{energydis} can be regarded as minimizing the energy dissipation of the resistive network interpreted by $\mG(\mV,\mE,\bm{W})$ when the total current injection at $\mV_a$ is unity (see \eqref{energydis1}), the total current ejection at $\mV_b$ is unity (see \eqref{energydis2}) and the other nodes have zero current injections (see \eqref{energydis3}).
By Theorem \ref{thmenergyVdiff}, the minimum energy dissipation in this case is $\reff(\Lg, \mV_a, \mV_b)$.
From this perspective, the resistive network can be equivalent to a resistance $\reff(\Lg, \mV_a, \mV_b)$ connecting two nodes that represent the aggregation of $\mV_a$ and $\mV_b$, respectively, see Fig. \ref{figgeffenergy} for illustration.
\end{remark}

\begin{figure}[!h]
  \centering
  \includegraphics[width=3.3in]{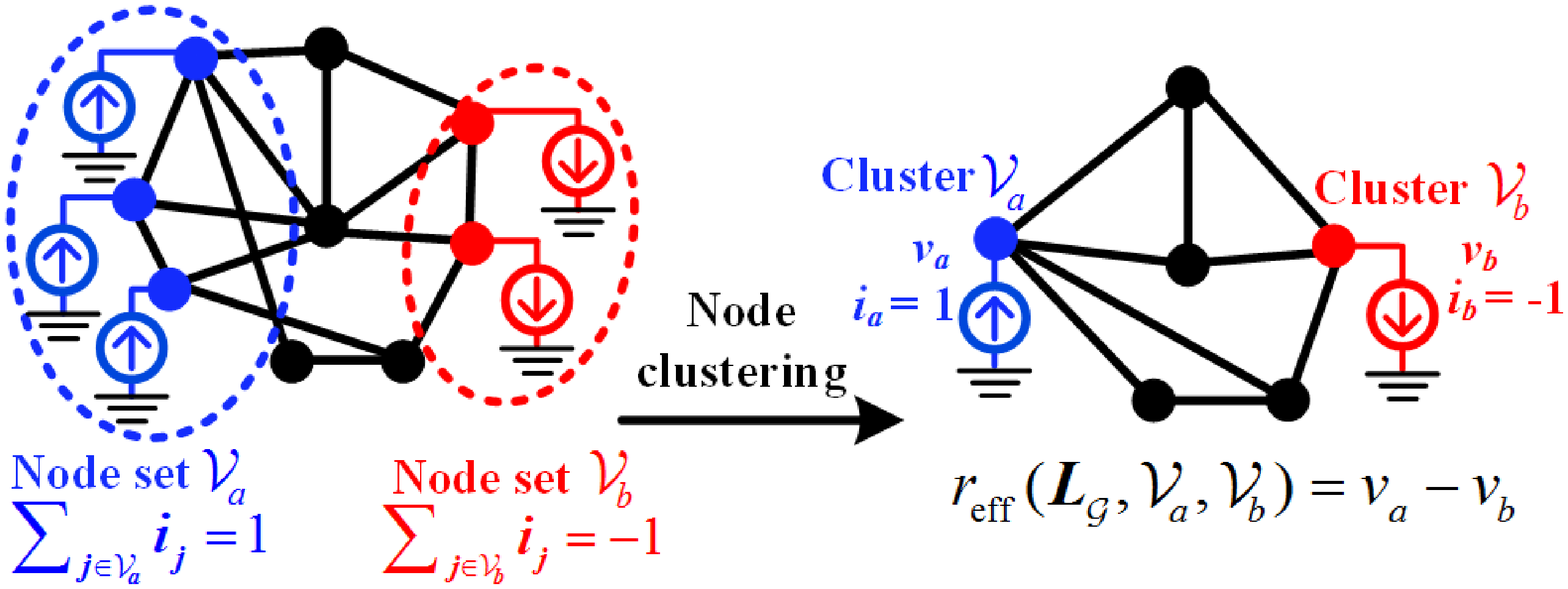}
  \caption{Interpretation of $\reff(\Lg, \mV_a, \mV_b)$ by potential difference.}
  \label{figgeffpotential}
\end{figure}

\begin{remark}
Equation \eqref{Vdiff} interprets $\reff(\Lg, \mV_a, \mV_b)$ by potential difference.
The matrix $\bm{\mP}_{ab}^T\Lg\bm{\mP}_{ab}$ is the Laplacian matrix of the graph where $\mV_a, \mV_b$ are clustered into two nodes (say node 1 and node 2) with the external connections being unchanged.
When one cluster has unit current injection and the other has unit current ejection, the potential difference between the two clusters is given by $(\bm{e}_1-\bm{e}_2)^T(\bm{\mP}_{ab}^T\Lg\bm{\mP}_{ab})^{\dag}(\bm{e}_1-\bm{e}_2)$, which is equal to $\reff(\Lg, \mV_a, \mV_b)$ by Theorem \ref{thmenergyVdiff}, see Fig. \ref{figgeffpotential} for illustration.
\end{remark}

\begin{remark}
Theorem \ref{thmenergyVdiff} leads to an interesting observation on Kron reduction and node clustering.
Let $\bm{L}_{cc}^{\mP}$ be the submatrix of $\bm{\mP}_{ab}^T\Lg\bm{\mP}_{ab}$ whose rows and columns are indexed by $\mV_c$. Simple circuit calculation on the clustered graph in Fig. \ref{figgeffpotential} gives that $\geff(\Lg, \mV_a, \mV_b)=\bm{e}_1^T[(\bm{\mP}_{ab}^T\Lg\bm{\mP}_{ab})/ \bm{L}_{cc}^{\mP}]\bm{e}_1$, which corresponds to Kron reduction on $\mV_c$ after node clustering on $\mV_a, \mV_b$.
On the other hand, the original definition of $\geff(\Lg, \mV_a, \mV_b)$ in \eqref{equgeff} can be regarded as the node clustering on $\mV_a, \mV_b$ after implementing Kron reduction on $\mV_c$.
It implies that the operation orders of Kron reduction on $\mV_c$ and node clustering on $\mV_a, \mV_b$ are exchangeable.
\end{remark}

\begin{figure}[!h]
  \centering
  \includegraphics[width=1.0in]{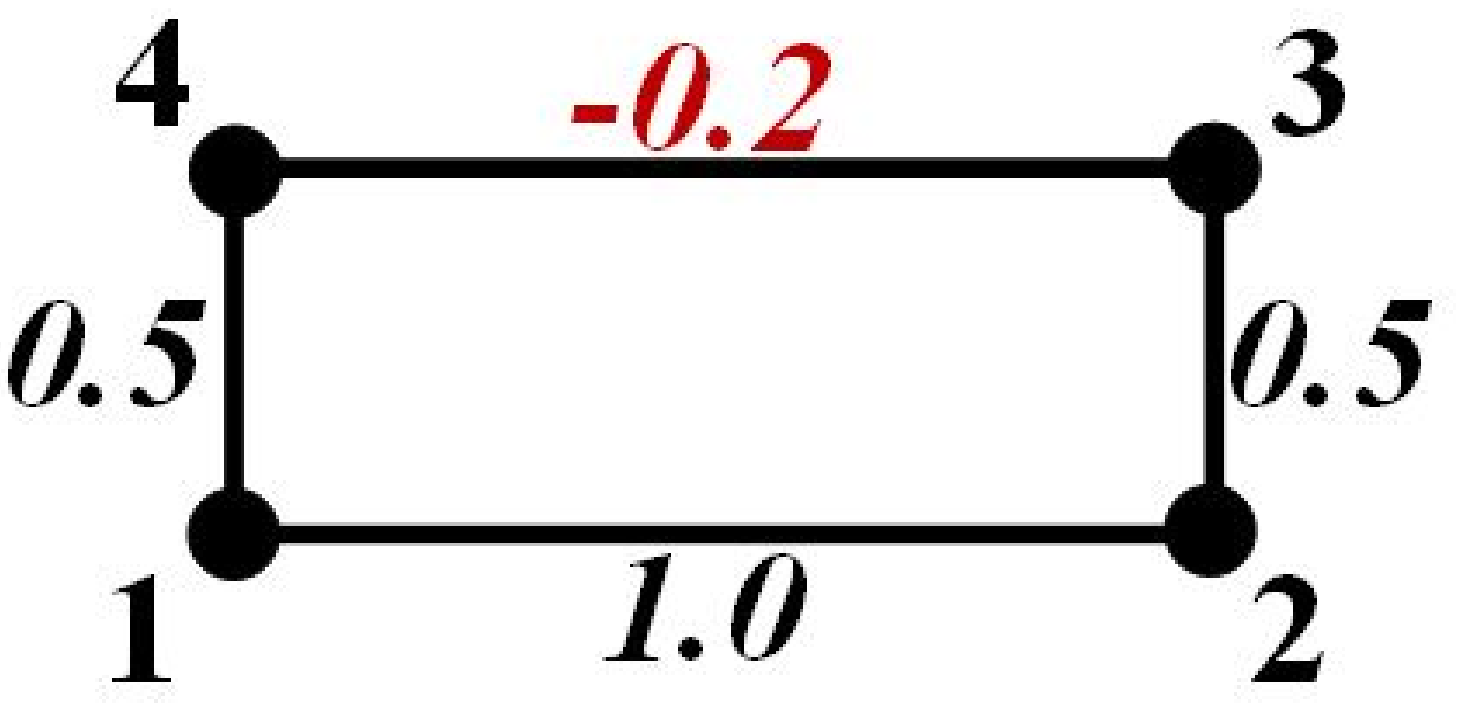}
  \caption{An example graph with negative weighted edges.}
  \label{figgeffmerit}
\end{figure}

\begin{remark}
The extended definitions outperform the traditional ones when applied to signed graphs.
To see this, consider a resistive network interpreted by the graph in Fig. \ref{figgeffmerit}, where the italic numbers denote line conductances.
The negative line conductance can be resulted from the port characteristics of some active circuit devices \cite{chua1987linear}.
The effective conductance between node 1 and node 2 can be given by simple circuit calculation $w_{12}+(w_{14}^{-1}+w_{34}^{-1}+w_{23}^{-1})^{-1}=0$. It implies that the negative line conductance causes ``electrical disconnection'' between node 1 and node 2.
Meanwhile, we have $\geff(\Lg, \{1\},\{2\})=0$ and $\mC_{12}=1.352$, which means the extended definition provides the correct answer but the traditional definition does not. The traditional definition fails due to $\Lg$ having two zero eigenvalues in this case and Lemma \ref{lemLdag} no longer holds.
\end{remark}

\indent
The above remarks conclude that $\reff(\Lg, \mV_a, \mV_b)$ and $\geff(\Lg, \mV_a, \mV_b)$ are a proper characterization of the overall coupling between two sets of nodes, which generalize the traditional ones with similar circuit interpretations from the perspective of current flow, energy dissipation and potential difference.
A smaller positive $\reff(\Lg, \mV_a, \mV_b)$ (or equivalently, larger positive $\geff(\Lg, \mV_a, \mV_b)$) indicates a stronger coupling strength.
Also note that $\geff(\Lg, \mV_a, \mV_b)$ takes the same physical dimension as the edge weight, which can be regarded as an ``equivalent weight''. Recall that a negative weighted edge is referred to as an antagonistic interaction between the two ternimal nodes in the consensus protocols \cite{altafini2013consensus}.
Borrowing this terminology, a negative $\reff(\Lg, \mV_a, \mV_b)$ or $\geff(\Lg, \mV_a, \mV_b)$ indicates that the two sets of nodes present ``antagonism'' to each other, and furthermore, a smaller negative $\geff(\Lg, \mV_a, \mV_b)$ indicates a more serious antagonism.

\subsection{Monotonicity and convexity}\label{secmonotone}
The traditional effective resistance and effective conductance have properties of monotonicity and convexity \cite{ghosh2008minimizing}.
The following results reveal that $\reff(\Lg, \mV_a, \mV_b)$ and $\geff(\Lg, \mV_a, \mV_b)$ preserve these properties.

\begin{theorem}[Monotonicity w.r.t. Node Sets]\label{thmmonotonenode}
   If $\Lg$ is PSD with only one zero eigenvalue, then $0<\reff(\Lg, \mV_a, \mV_b)\leq \reff(\Lg, \mV_a^{*}, \mV_b^{*})$ and $\geff(\Lg, \mV_a, \mV_b)\geq \geff(\Lg, \mV_a^{*}, \mV_b^{*})>0$ for any $\mV_a, \mV_b, \mV_a^{*}, \mV_b^{*}$ such that $\mV_a, \mV_b\subset \mV$, $\mV_a\cap \mV_b=\phi$, $\mV_a^{*}\subseteq \mV_a$ and $\mV_b^{*}\subseteq \mV_b$.
\end{theorem}

\begin{IEEEproof}
    By Theorem \ref{thmenergyVdiff}, $\reff(\Lg, \mV_a^{*}, \mV_b^{*})$ and $\reff(\Lg, \mV_a, \mV_b)$ equal to the optimal values of the corresponding optimization problems \eqref{energydis}.
    When $\mV_a^{*}\subseteq \mV_a$ and $\mV_b^{*}\subseteq \mV_b$, the feasible set of problem \eqref{energydis} with respect to $\reff(\Lg, \mV_a^{*}, \mV_b^{*})$ is a subset of that with respect to $\reff(\Lg, \mV_a, \mV_b)$. So we have $\reff(\Lg, \mV_a, \mV_b)\leq \reff(\Lg, \mV_a^{*}, \mV_b^{*})>0$ and $\geff(\Lg, \mV_a, \mV_b)\geq \geff(\Lg, \mV_a^{*}, \mV_b^{*})>0$.
\end{IEEEproof}

\begin{theorem}[Monotonicity w.r.t. Edge Weights]\label{thmmonotone}
   If $\Lg$ is PSD with only one zero eigenvalue, then for any two non-empty disjoint sets $\mV_a, \mV_b\subset \mV$, $\geff(\Lg, \mV_a, \mV_b)$ is non-decreasing and $\reff(\Lg, \mV_a, \mV_b)$ is non-increasing after adding an edge $(i,j)$ with weight $\Delta w_{ij}>0$ to the graph $\mG$.
\end{theorem}

\begin{IEEEproof}
    The proof is given in Appendix B.
\end{IEEEproof}

\begin{theorem}[Convexity]\label{thmconvexity}
   For any two non-empty disjoint sets $\mV_a, \mV_b\subset \mV$, $\geff(\Lg, \mV_a, \mV_b)$ is a concave function of $\bm{W}$ and $\reff(\Lg, \mV_a, \mV_b)$ is a convex function of $\bm{W}$ if $\Lg$ is PSD with only one zero eigenvalue.
\end{theorem}

\begin{IEEEproof}
   We first prove the concavity of $\geff(\Lg, \mV_a, \mV_b)$.
   Let $\bm{E}_a\in \mbR^{a\times l}$, $\bm{E}_b\in \mbR^{b\times l}$, $\bm{E}_c\in \mbR^{c\times l}$ denote the submatrices of the incidence matrix $\bm{E}$ whose rows are indexed by $\mV_a, \mV_b, \mV_c$, respectively. Then, $\Lg$ can be rewritten as
   \begin{equation*}
   \small
   \begin{split}
     \Lg =
     \begin{bmatrix}
     \bm{L}_{aa}   &   \bm{L}_{ab}   &   \bm{L}_{ac}     \\
     \bm{L}_{ab}^T   &   \bm{L}_{bb}   &   \bm{L}_{bc}    \\
      \bm{L}_{ac}^T  &   \bm{L}_{bc}^T   &   \bm{L}_{cc}    \\
     \end{bmatrix} =
     \begin{bmatrix}
      \bm{E}_a\bm{W}\bm{E}_a^T   &   \bm{E}_a\bm{W}\bm{E}_b^T   &   \bm{E}_a\bm{W}\bm{E}_c^T     \\
      \bm{E}_b\bm{W}\bm{E}_a^T   &   \bm{E}_b\bm{W}\bm{E}_b^T   &   \bm{E}_b\bm{W}\bm{E}_c^T    \\
      \bm{E}_c\bm{W}\bm{E}_a^T  &   \bm{E}_c\bm{W}\bm{E}_b^T   &   \bm{E}_c\bm{W}\bm{E}_c^T    \\
     \end{bmatrix}
   \end{split}
   \end{equation*}
   and by \eqref{Lcc3aa} we have
   \begin{equation*}
   \begin{split}
      \geff(\Lg, \mV_a, \mV_b) =& \bm{1}_a^T\bm{E}_a\bm{W}\bm{E}_a^T\bm{1}_a - \\ &\bm{1}_a^T\bm{E}_a\bm{W}\bm{E}_c^T(\bm{E}_c\bm{W}\bm{E}_c^T)^{-1}\bm{E}_c\bm{W}\bm{E}_a^T\bm{1}_a.
   \end{split}
   \end{equation*}
   Define the hypograph \cite{boyd2004convex} of $\geff(\Lg, \mV_a, \mV_b)$ as
   \begin{equation*}
   \begin{split}
   \textup{hypo}\{\geff(\Lg, \mV_a, \mV_b)\} = \{(\bm{W}, t)|~t\leq \geff(\Lg(\bm{W}), \mV_a, \mV_b)\}.
   \end{split}
   \end{equation*}
   Since $\Lg$ is PSD with only one zero eigenvalue, we have $\bm{L}_{cc}=\bm{E}_c\bm{W}\bm{E}_c^T\succ 0$.
   Then the term $\geff(\Lg, \mV_a, \mV_b)-t$ can be regarded as a Schur complement $\bm{H}/ \bm{E}_c\bm{W}\bm{E}_c^T$, where
   \begin{equation}\label{geffschur}
   \begin{split}
      \bm{H} =
      \begin{bmatrix}
         \bm{1}_a^T\bm{E}_a\bm{W}\bm{E}_a^T\bm{1}_a-t & \bm{1}_a^T\bm{E}_a\bm{W}\bm{E}_c^T \\
         \bm{E}_c\bm{W}\bm{E}_a^T\bm{1}_a & \bm{E}_c\bm{W}\bm{E}_c^T \\
      \end{bmatrix}.
   \end{split}
   \end{equation}
   Thus, by Lemma \ref{lemschur}, the hypograph is equivalent to $\textup{hypo}\{\geff(\Lg, \mV_a, \mV_b)\} = \{(\bm{W}, t)|~\bm{H}\succeq 0 \}$, which is a convex set with respect to $(\bm{W}, t)$. So $\geff(\Lg, \mV_a, \mV_b)$ is a concave function of $\bm{W}$ according to \cite{boyd2004convex}.

   \indent
   We have proved that $\geff(\Lg, \mV_a, \mV_b)$ is a concave function of $\bm{W}$ and $\geff(\Lg, \mV_a, \mV_b)>0$ (see Theorem \ref{thmenergyVdiff}) if $\Lg$ is PSD with only one zero eigenvalue.
   In addition, it is trivial that $h(x)=1/x$ is convex and non-increasing for $x>0$.
   Then, by the convexity condition for composite functions \cite{boyd2004convex}, $\reff(\Lg, \mV_a, \mV_b)=h(\geff(\Lg, \mV_a, \mV_b))$ is convex of $\bm{W}$.
\end{IEEEproof}

\indent
Generally, the results obtained in this subsection tell that $\reff(\Lg, \mV_a, \mV_b)$ (or $\geff(\Lg, \mV_a, \mV_b)$) is decreasing and convex (or increasing and concave) with respect to edge weights in case that $\Lg$ is PSD with only one zero eigenvalue. These properties will be leveraged to establish further results in the following sections.

\section{Interpreting graph Laplacian definiteness by effective resistance}\label{seclapPSD}
As mentioned in the introduction section, graph Laplacian being PSD with only one zero eigenvalue is a key property required by many problems over networks, such as reaching a consensus or synchronization \cite{zelazo2014definiteness, dorfler2018electrical}.
However, this property may not be held by signed graphs, which depends on the locations and weights of the negative weighted edges.
In this section, we will show that the sign of effective resistance or effective conductance provides rich information for the definiteness of graph Laplacian.
Since $\reff(\Lg, \mV_a, \mV_b)$ and $\geff(\Lg, \mV_a, \mV_b)$ are reciprocal to each other, the following results will be stated in terms of $\geff(\Lg, \mV_a, \mV_b)$ for simplicity. In addition, unless otherwise specified, the terms ``effective resistance'' and ``effective conductance'' refer to the extended definitions henceforth.

\indent
Given a graph $\mG(\mV, \mE, \bm{W})$, let $\mG_-(\mV, \mE_-, \bm{W}_-)$ and $\mG_+(\mV, \mE_+, \bm{W}_+)$ be its negative subgraph and positive subgraph, where $\mE_-$ and $\mE_+=\mE\backslash \mE_-$ denote the set of edges with negative weights and positive weights, respectively.
Denote $\mV_-=\{i|~\exists (i,j)\in\mE:~w_{ij}<0\}$ as the set of nodes that are adjacent to negative weighted edges.
Also we propose the following definition.

\begin{definition}[Sequential Inclusion of Set]
   For a set $\mV$ with cardinality $n$, define the group of subsets $\{\mV^k\}$, $k=1,2,...,n$ as a sequential inclusion of $\mV$ if $\mV^1\subset\mV^2\subset\cdots\subset\mV^n$, $\mV^n=\mV$ and $|\mV^{k+1}|-|\mV^{k}|=1$, $k=1,...,n-1$.
\end{definition}

\indent
For instance, $\mV^1=\{2\}$, $\mV^2=\{2,3\}$ and $\mV^3=\{2,3,1\}$ form a sequential inclusion of $\mV=\{1,2,3\}$.

\indent
With these notations, we have the following results on graph Laplacian definiteness.

\begin{lemma}\label{lemPSD}
   $\Lg$ is PSD with only one zero eigenvalue if there exists a sequential inclusion of $\mV$, say $\{\mV^k\}$, such that
   \begin{equation}\label{geffpos}
      \geff(\Lg, \mV^{k+1}\backslash\mV^k, \mV^k)>0,~k=1,2,...,n-1.
   \end{equation}
\end{lemma}

\begin{IEEEproof}
   The proof is given in Appendix C.
\end{IEEEproof}

\begin{theorem}[Laplacian Definiteness \& Effective Conductance]\label{thmPSD}
   Given a graph $\mG(\mV, \mE, \bm{W})$, the following statements are equivalent:
   \begin{description}
     \item[(a)] $\Lg$ is PSD with only one zero eigenvalue.
     \item[(b)] There exists a sequential inclusion of $\mV_-$, say $\{\mV_-^k\}$, such that $\geff(\Lg, \mV_-^{k+1}\backslash\mV_-^k, \mV_-^k)>0$, $k=1,2,...,|\mV_-|-1$.
     \item[(c)] There exists a set $\mV_o\subset \mV$ and a sequential inclusion of $\mV_r=\mV\backslash\mV_o$, say $\{\mV_r^k\}$, such that $\bm{L}_{oo}$ is positive definite and $\geff(\Lg, \mV_r^{k+1}\backslash\mV_r^k, \mV_r^k)>0$, $k=1,2,...,|\mV_r|-1$.
   \end{description}
\end{theorem}

\begin{IEEEproof}
    First, it is trivial to show that (a)$\Rightarrow$(b) and (a)$\Rightarrow$(c) follow from Theorem \ref{thmenergyVdiff}.

    \indent
    Second, we prove (b)$\Rightarrow$(a). Let $\bm{L}_{++}$ denote the submatrix of $\Lg$ whose rows and columns are indexed by $\mV_+=\mV\backslash\mV_-$. It is trivial that $\bm{L}_{++}$ is positive definite.
    Thus, by Lemma \ref{lemschur} we have $i_-(\Lg)=i_-(\Lg/ \bm{L}_{++})$ and $i_0(\Lg)=i_0(\Lg/ \bm{L}_{++})$. Let $\mG_r(\mV_r, \mE_r, \bm{W}_r)$, $\mV_r=\mV_-$ be the graph whose Laplacian matrix is $\Lg/ \bm{L}_{++}$.
    By Lemma \ref{lemquotient} we have $\geff(\Lgr, \mV_-^{k+1}\backslash\mV_-^k, \mV_-^k)=\geff(\Lg, \mV_-^{k+1}\backslash\mV_-^k, \mV_-^k)>0$ for the sequential inclusion of $\mV_-$.
    Then, it follows from Lemma \ref{lemPSD} that $\Lg/ \bm{L}_{++}$ is PSD with only one zero eigenvalue, and hence $\Lg$ is also PSD with only one zero eigenvalue.

    \indent
    Third, we prove (c)$\Rightarrow$(a). Since $\bm{L}_{oo}$ is positive definite, it follows from Lemma \ref{lemschur} that $i_-(\Lg)=i_-(\Lg/ \bm{L}_{oo})$ and $i_0(\Lg)=i_0(\Lg/ \bm{L}_{oo})$. Then, we can conclude that $\Lg$ is PSD with only one zero eigenvalue by using the same logic as in the proof of (b)$\Rightarrow$(a).
\end{IEEEproof}

\begin{corollary}\label{corPSD}
   $\Lg$ is PSD with only one zero eigenvalue if and only if $\geff(\Lg, \mV_a, \mV_b)> 0$ for any two non-empty disjoint sets $\mV_a, \mV_b\subset \mV$.
\end{corollary}

\begin{IEEEproof}
   The necessity follows from Theorem \ref{thmenergyVdiff}, and the sufficiency follows from Theorem \ref{thmPSD} by letting $\mV_o=\phi$ in its statement (c).
\end{IEEEproof}

\begin{theorem}[Laplacian Negative Eigenvalues \& Negative Weighted Edges]\label{thmnegative1}
    $\Lg$ has at most $|\mE_-|$ negative eigenvalues.
\end{theorem}

\begin{IEEEproof}
    For any graph $\mG_0(\mV,\mE_0,\bm{W}_0)$ (not necessarily be identical to $\mG$), let us consider graph $\mG_*$ that is obtained by adding an edge $(i,j)$ to $\mG_0$.
    Let $\mV_a=\{i\}$, $\mV_b=\{j\}$ and $\mV_c=\mV\backslash\{i,j\}$.
    Since $\Lgs/\bm{L}_{cc}$ is a two-by-two matrix in this case, by Lemma \ref{lemLcc} we have $\Lgs/\bm{L}_{cc}=\geff(\Lgs, \{i\}, \{j\})\cdot
   \begin{bmatrix}
   1 & -1 \\
   -1 & 1 \\
   \end{bmatrix}$, where $\geff(\Lgs, \{i\}, \{j\})$ is the effective conductance between node $i$ and node $j$ in graph $\mG_*$, and $\bm{L}_{cc}$ is a submatrix of $\Lgs$ as well as $\Lgz$.
   Further, by Lemma \ref{lemschur} we have $i_-(\Lgs) = i_-(\bm{L}_{cc}) + i_-(\Lgs/\bm{L}_{cc})$ and hence
   \begin{equation}\label{lapinertia}
   \begin{split}
       i_-(\Lgs) \leq i_-(\Lgz) + i_-(\geff(\Lgs, \{i\}, \{j\}))
   \end{split}
   \end{equation}
   which implies that the number of negative eigenvalues of graph Laplacian is increased by at most one after adding a new edge.

   \indent
   We now turn to the graph $\mG(\mV,\mE,\bm{W})$.
   Note that $\mG$ can be obtained by adding the set of edges $\mE_-$ to its positive subgraph $\mG_+(\mV, \mE_+, \bm{W}_+)$, and the Laplacian matrix of $\mG_+$ has no negative eigenvalues.
   Thus, as inferred from \eqref{lapinertia}, $\Lg$ has at most $|\mE_-|$ negative eigenvalues.
\end{IEEEproof}

\begin{remark}
Theorem \ref{thmPSD} gives two necessary and sufficient conditions for $\Lg$ being PSD with only one zero eigenvalue.
It generalizes the results in \cite{zelazo2014definiteness, chen2016definiteness, zelazo2017robustness} that apply to those graphs with a special location distribution of negative weighted edges.
Note that the sequential inclusion describes an expansion from a single node to the whole concerned set ($\mV_-^{k+1}\backslash\mV_-^k$ refers to the node to be included in one step). The condition $\geff(\Lg, \mV_-^{k+1}\backslash\mV_-^k, \mV_-^k)>0$ indicates that no antagonisms have occurred during such an expansion.
From this viewpoint, statement (b) of Theorem \ref{thmPSD} gives a twofold description for $\Lg$ being PSD with only one zero eigenvalue.
First, there is a set of nodes $\mV_+=\mV\backslash\mV_-$ that are not impacted by negative weighted edges and the submatrix $\bm{L}_{++}$ is positive definite. Second, the remaining set of nodes $\mV_-$, which are impacted by negative weighted edges, can be formed by a sequential inclusion without antagonisms.
Statement (c) of Theorem \ref{thmPSD} can be regarded as a further generalization of statement (b).
\end{remark}

\begin{remark}
Corollary \ref{corPSD} provides a neat form that elaborates the mechanism of graph Laplacian definiteness.
The Laplacian being PSD with only one zero eigenvalue can be interpreted as positive couplings everywhere in the graph, i.e., the negative weighted edges do not induce antagonism between any two sets of nodes.
\end{remark}

\begin{remark}
It is previously proved in \cite{pan2016laplacian} that the number of Laplacian negative eigenvalues of a cycle graph is upper bounded by $|\mE_-|$.
Theorem \ref{thmnegative1} extends this result to generic graphs based on the properties of effective conductance.
\end{remark}

\begin{remark}
Theorem \ref{thmPSD} and Corollary \ref{corPSD} also motivate an algorithm below for checking if $\Lg$ is PSD with only one zero eigenvalue.
\begin{enumerate}[Step 1:]
  \item Initialization. Arbitrarily choose a node $i\in\mV_-$, set $\mV_b=\{i\}$ and update $\mV_-\leftarrow \mV_-\backslash \{i\}$.
  \item Arbitrarily choose a node $j\in\mV_-$, and update $\mV_-\leftarrow \mV_-\backslash \{j\}$.
  \item If $\geff(\Lg, \{j\}, \mV_b)>0$, update $\mV_b\leftarrow \mV_b\cup \{j\}$.
  Otherwise stop the algorithm, $\Lg$ is not PSD with only one zero eigenvalue.
  \item Stop if $\mV_-= \phi$, $\Lg$ is PSD with only one zero eigenvalue.
  Otherwise go back to Step 2.
\end{enumerate}
This algorithm is efficient as $\mV_-$ usually has small cardinality in practical problems.
\end{remark}

\section{Interpreting power network stability by effective resistance}\label{secstability}
\subsection{Indication of small-disturbance angle stability}\label{secsds}
Small-disturbance angle stability refers to the ability of a power system to maintain synchronism under small disturbances \cite{kundur2004definition}, which is a fundamental requirement for power system operation. Small-disturbance stability analysis also provides critical information for the geometry of stability boundary, e.g., the type of an unstable equilibrium point (UEP).
The existing results on angle stability are much oriented to node dynamics, particularly generator dynamics.
In the following, we will show that the effective conductance is a powerful stability indicator and leads to a network-based interpretation of angle stability.

\indent
Consider a structure-preserving power network\footnote{The power networks considered here are different from the electrical networks used for the illustration of effective resistance in the previous sections. The former one refers to high-voltage AC transmission networks with inductive lines. The latter one refers to DC networks with resistive lines.} where the set of buses\footnote{For terminology convention in the respective disciplines, we interchangeably use ``buses, lines'' for power networks and ``nodes, edges'' for graphs.} $\mV$ are interconnected via the set of transmission lines $\mE$.
The generators are modeled as internal buses (voltage sources) connecting to the terminal buses via transient reactances, and the loads are frequency dependent. These fictitious internal buses and transient reactances are included into $\mV$ and $\mE$, respectively.
The transmission lines are assumed purely inductive as the resistances of physical lines are negligible for high-voltage transmission networks.
Then, the system dynamics can be described by \cite{bergen1981structure}
\begin{subequations}\label{tsa0}
\begin{align}
   M_i\ddot{\theta}_i + D_i\dot{\theta}_i&= P_i - \sum\limits_{j\in\mN_i} V_iV_jB_{ij}\sin\theta_{ij},~i\in\mV_G \label{tsa0g} \\
   D_i\dot{\theta}_i &= P_i - \sum\limits_{j\in\mN_i} V_iV_jB_{ij}\sin\theta_{ij} ,~i\in\mV_L  \label{tsa0l}
\end{align}
\end{subequations}
where $\mV_G=\{1,...,g\}$ and $\mV_L=\mV\backslash\mV_G=\{g+1,...,n\}$ denote the set of generator internal buses and the remaining buses, respectively (let $|\mV_G|=g$, $|\mV_L|=d$);
$\mN_i=\{j|(i,j)\in \mE\}$ denotes the set of neighboring buses of bus $i$;
$\theta_i, V_i$ denote the rotor angle and voltage magnitude of bus $i$, respectively; $\theta_{ij}$ is defined by $\theta_{ij}=\theta_i-\theta_j$; the susceptance of line $(i,j)$ is equal to $-B_{ij}$;
$P_i$ denotes the active power injection at bus $i$; $M_i>0$ denotes the generator inertia; and $D_i>0$ denotes the damping coefficient for bus $i\in\mV_G$ or load frequency sensitivity for $i\in \mV_L$.
We assume $V_i$ keeps constant and $\sum_{i\in \mV} P_i=0$, the latter of which gives $\dot{\theta}_i=0$, $\forall i\in \mV$ at an equilibrium point. The above modeling and assumptions are common in the study of angle stability (e.g., see \cite{song2017cutset} for a detailed discussion). Also, we introduce the concepts below.

\begin{definition}[Active Power Flow Graph \& Critical Line \cite{song2017networkbased}]\label{defAPFG}
   The graph $\mG(\mV, \mE, \bm{W}(\bm{\theta}))$ is called the active power flow graph with the set of buses $\mV$, the set of lines $\mE$ and the diagonal matrix of line weights $\bm{W}(\bm{\theta})=\diag{V_iV_jB_{ij}\cos\theta_{ij}}\in \mbR^{l\times l}$, $\forall (i,j)\in \mE$. A line $(i,j)$ in the active power flow graph $\mG(\mV, \mE, \bm{W}(\bm{\theta}))$ is called a critical line if $\frac{\pi}{2}<|\theta_i-\theta_j| \mod 2\pi < \frac{3\pi}{2}$.
\end{definition}

\indent
The active power flow graph has the same topology as the physical power network, and the edge weights are characterized by line susceptances and angle differences.
The set of critical lines is identical to the set of negative weighted edges in the active power flow graph, denoted by $\mE_-$.
The critical lines are likely to occur at those equilibrium points with large angle differences across some lines, which refer to heavy load level or abnormal operating scenario.
These concepts will be used to connect the obtained theorems on graph Laplacian definiteness to small-disturbance angle stability.

\indent
With the active power graph, we obtain the small-disturbance model below by linearizing \eqref{tsa0} around an equilibrium point
\begin{equation}\label{sda}
\begin{split}
   \begin{bmatrix}
     \bm{M}\Delta\ddot{\bm{\theta}}_G \\
     \bm{0}_d  \\
   \end{bmatrix} +
   \begin{bmatrix}
     \bm{D}_G\Delta\dot{\bm{\theta}}_G \\
     \bm{D}_L\Delta\dot{\bm{\theta}}_L \\
   \end{bmatrix} +
   \Lg(\bm{\theta}^e)\Delta\bm{\theta}=\bm{0}_n
\end{split}
\end{equation}
where $\bm{\theta}_G=[\theta_i]\in \mbR^g$, $\forall i\in \mV_G$;
$\bm{\theta}_L=[\theta_i]\in \mbR^d$, $\forall i\in \mV_L$;
$\bm{\theta}=\begin{bmatrix}
               \bm{\theta}_G^T & \bm{\theta}_L^T \\
             \end{bmatrix}^T\in\mbR^n$;
$\bm{D}_G=\diag{D_i}\in \mbR^{g\times g}$, $\forall i\in \mV_G$;
$\bm{D}_L=\diag{D_i}\in \mbR^{d\times d}$, $\forall i\in \mV_L$;
the variable with the superscript ``$e$'' denotes the value at the equilibrium point;
and $\Lg(\bm{\theta}^e)$ denotes the Laplacian matrix of $\mG(\mV, \mE, \bm{W}(\bm{\theta}^e))$.
Further, by taking bus $n$ as the angle reference, the small-disturbance model \eqref{sda} can be re-expressed as the following state space form \cite{bergen1981structure}
\begin{equation}\label{sdastate}
\begin{split}
    \begin{bmatrix}
      \Delta\dot{\bm{\alpha}} \\
      \Delta\dot{\bm{\omega}}_G \\
    \end{bmatrix}
    =&
    \begin{bmatrix}
      -\bm{T}_L\bm{D}_L^{-1}\bm{T}_L^{T}\bm{F}(\bm{\theta}^e) & \bm{T}_G \\
      -\bm{M}^{-1}\bm{T}_G^{T}\bm{F}(\bm{\theta}^e) & -\bm{M}^{-1}\bm{D}_G \\
    \end{bmatrix}
    \begin{bmatrix}
      \Delta\bm{\alpha} \\
      \Delta\bm{\omega}_G \\
    \end{bmatrix}
    \\=&\bm{J}_{dyn}(\bm{\theta}^e)
    \begin{bmatrix}
      \Delta\bm{\alpha} \\
      \Delta\bm{\omega}_G \\
    \end{bmatrix}
\end{split}
\end{equation}
where $\bm{T}_G\in \mbR^{(n-1)\times g}$ and $\bm{T}_L\in \mbR^{(n-1)\times d}$ are submatrices of
$\bm{T}=
\begin{bmatrix}
       \bm{I}_{n-1} & -\bm{1}_{n-1} \\
\end{bmatrix}\in \mbR^{(n-1)\times n}$ whose columns are indexed by $\mV_G$ and $\mV_L$, respectively;
$\bm{\alpha}=\bm{T}\bm{\theta}\in \mbR^{n-1}$ denotes the vector of relative angles with respect to bus $n$;
$\bm{\omega}_G=\dot{\bm{\theta}}_G\in \mbR^g$ denotes the vector of generator rotor speeds;
and $\bm{F}(\bm{\theta}^e)\in \mbR^{(n-1)\times(n-1)}$ is the Laplacian matrix $\Lg(\bm{\theta}^e)$ with the last row and column deleted.
Given an equilibrium point of \eqref{sdastate}, say $(\bm{\alpha}^{e}, \bm{0}_g)$, and the corresponding angle vector
$\bm{\theta}^e=\begin{bmatrix}
                               (\bm{\alpha}^e)^T & 0 \\
                        \end{bmatrix}^T$,
we have the following concepts.

\begin{definition}[Hyperbolicity, Stability and Type of Equilibria \cite{chiang1987foundations}]\label{defSEP}
    An equilibrium point $(\bm{\alpha}^e, \bm{0}_g)$ is hyperbolic if $i_0(\bm{J}_{dyn}(\bm{\theta}^e))=0$.
    For a hyperbolic equilibrium point, it is locally asymptotically stable if $i_+(\bm{J}_{dyn}(\bm{\theta}^e))=0$, otherwise, it is a type-$m$ UEP where $m=i_+(\bm{J}_{dyn}(\bm{\theta}^e))$.
\end{definition}

\indent
The type of UEP refers to the number of unstable manifolds around a UEP, which is an important geometric description indicating possible directions for system states to escape the stability region \cite{chiang1987foundations}.

\indent
We come to the following theorems that characterize hyperbolicity, stability and type of the equilibrium point in terms of $\geff(\Lg(\bm{\theta}^e), \mV_a, \mV_b)$.
Also note that $\geff(\Lg(\bm{\theta}^e), \mV_a, \mV_b)$ actually describes the effective susceptance (i.e., electrical coupling strength) in the context of power networks as the edge weights in the active power flow graph are in terms of $B_{ij}$ and $\theta_{ij}^e$.
Nevertheless, we still refer to it as effective conductance to avoid confusion with the original definition.

\begin{theorem}[Stability \& Graph Laplacian]\label{thmSEPandL}
   The equilibrium point $(\bm{\alpha}^e, \bm{0}_g)$ is hyperbolic and locally asymptotically stable if and only if $\Lg(\bm{\theta}^e)$ is PSD with only one zero eigenvalue.
\end{theorem}

\begin{IEEEproof}
   First we claim $i_-(\bm{F})=i_-(\Lg)$ and $i_0(\Lg)=i_0(\bm{F})+1$, which can be proved by applying Sylvester's law of inertia \cite{horn2012matrix} to $\Lg=\bm{S}\bm{K}\bm{S}^{T}$, where
   $\bm{K} = \begin{bmatrix}
                         \bm{F} & \bm{0}_{n-1} \\
                               \bm{0}_{n-1}^T & 0 \\
                     \end{bmatrix},~
   \bm{S} = \begin{bmatrix}
                          \bm{I}_{n-1} & \bm{0}_{n-1} \\
                               -\bm{1}_{n-1}^{T} & 1 \\
                     \end{bmatrix}$.
   The term ``$(\bm{\theta}^e)$'' is dropped in the concerned matrices for simplicity.

   \indent
   We now prove the sufficiency.
   If $\Lg$ has only one zero eigenvalue, then $\bm{F}$ is nonsingular due to $i_0(\Lg)=i_0(\bm{F})+1$, and thus the equilibrium point is hyperbolic by [\citen{chiang1987foundations}, Proposition 3].
   Further, it follows from [\citen{chiang1987foundations}, Theorem 7] that the hyperbolic equilibrium point satisfies $i_+(\bm{J}_{dyn})=i_-(\bm{F})$.
   Since $i_-(\bm{F})=i_-(\Lg)$ and $i_-(\Lg)=0$ ($\Lg$ is PSD), we have $i_+(\bm{J}_{dyn})=0$ and the equilibrium point is locally asymptotically stable.

   \indent
   Necessity. If the equilibrium point is hyperbolic, i.e., $i_0(\bm{J}_{dyn})=0$, then $\bm{F}$ is nonsingular by the relation between $\bm{J}_{dyn}$ and $\bm{F}$ in \eqref{sdastate}. So $\Lg$ has only one zero eigenvalue due to $i_0(\Lg)=i_0(\bm{F})+1$.
   Again, it follows from [\citen{chiang1987foundations}, Theorem 7] that $i_+(\bm{J}_{dyn})=i_-(\bm{F})$.
   Thus, if the equilibrium point is asymptotically stable, we have $i_-(\Lg)=i_-(\bm{F})=i_+(\bm{J}_{dyn})=0$ so that $\Lg$ is PSD.
\end{IEEEproof}

\begin{theorem}[Stability \& Effective Conductance]\label{thmSEPandgeff}
   For an equilibrium point $(\bm{\alpha}^e, \bm{0}_g)$ and the corresponding active power flow graph $\mG(\mV, \mE, \bm{W}(\bm{\theta}^e))$, the following statements are equivalent:
   \begin{description}
     \item[(a)] The equilibrium point is hyperbolic and locally asymptotically stable.
     \item[(b)] There exists a a sequential inclusion of $\mV_-$, say $\{\mV_-^k\}$, such that $\geff(\Lg(\bm{\theta}^e), \mV_-^{k+1}\backslash\mV_-^k, \mV_-^k)>0$, $k=1,2,...,|\mV_-|-1$.
     \item[(c)] $\geff(\Lg(\bm{\theta}^e), \mV_a, \mV_b)> 0$ for any two non-empty disjoint sets $\mV_a, \mV_b\subset \mV$.
   \end{description}
\end{theorem}

\begin{IEEEproof}
   It directly follows from Theorem \ref{thmPSD}, Corollary \ref{corPSD} and Theorem \ref{thmSEPandL}.
\end{IEEEproof}

\begin{theorem}[Type of UEP\& Critical Lines]\label{thmUEPandL}
   If the equilibrium point $(\bm{\alpha}^e, \bm{0}_g)$ is hyperbolic and unstable, then it is a type-$m$ UEP, where $m$ is less than or equal to the number of critical lines, i.e., $m\leq |\mE_-|$.
\end{theorem}

\begin{IEEEproof}
   According to the proof of Theorem \ref{thmSEPandL}, we have $i_+(\bm{J}_{dyn})=i_-(\Lg)$ if the equilibrium point is hyperbolic.
   Thus, we can complete the proof together with Theorem \ref{thmnegative1}.
\end{IEEEproof}

\begin{remark}
These theorems are direct consequences of those obtained in the previous section, which give a graph theoretic characterization for power network stability.
Theorem \ref{thmSEPandL} links small-disturbance angle stability to the definiteness of active power flow graph Laplacian.
Further, Theorem \ref{thmSEPandgeff} gives necessary and sufficient conditions for the hyperbolicity and asymptotic stability of an equilibrium point in terms of the effective conductance, which makes a generalization of our previous result [\citen{song2017networkbased}, Theorem 3].
Since the negative effective conductance in the active power flow graph can be interpreted as ``electrical antagonism'', Theorem \ref{thmSEPandgeff} indicates that angle instability is equivalent to the electrical antagonism between certain two sets of buses.
It will also be seen in case study that those pairs of bus sets with negative effective conductances tend to have angle separation when subjected to disturbances.
Moreover, Theorem \ref{thmUEPandL} implies that the number of unstable manifolds around a UEP is upper bounded by the number of critical lines.
It provides new ideas into linking the unstable manifolds to physical components in power systems.
\end{remark}

\subsection{Transient stability enhancement}\label{sectsa}
Transient stability refers to the ability of a power system to maintain synchronism under large disturbances \cite{kundur2004definition}.
In this subsection, we will design an effective conductance-based approach to enhancing the transient stability level of operating equilibrium point, which is the locally stable equilibrium used for normal operation with angle difference across each line being relatively small.
The analysis will be based on a common observation in power systems called generator coherency, i.e., (usually) two groups of generators have similar waveforms for their rotor angles after a fault \cite{wu1983identification}.

\indent
Consider the system runs at the operating equilibrium point and undergoes a fault which is cleared at time $t=0$.
Assume the fault is self-cleared so that the post-fault system has the same network topology as the pre-fault one, which is a typical scenario studied in the literature \cite{vu2016toward}.
The state variable of the post-fault system deviates from the operating equilibrium due to the fault, and we use the superscript ``$t$'' to denote those variables at time $t\geq 0$.
Let $\mV_a$ and $\mV_b=\mV_G\backslash \mV_a$ be the groups of coherent generators during the post-fault oscillation such that $|\theta_{ij}^t|<\varepsilon$, $\forall t\geq 0$, $\forall i,j\in \mV_a$ or $i,j\in \mV_b$ where $\varepsilon$ is a small positive number ($\mV_a,\mV_b$ can be identified by the method in \cite{wu1983identification}), and $\mV_c=\mV\backslash(\mV_a\cup\mV_b)$ be the set of remaining buses.
We also assume the effect of damping coefficients is excluded, i.e., $D_i=0$, $\forall i\in \mV$. This is a reasonable approximation for real situation as $D_i$ is practically small compared to $M_i$, which is often used in transient stability analysis \cite{kundur1994power}.
Then, the system is conservative \cite{bergen1981structure} and we focus on the interaction between the coherent generators.

\indent
With the above settings, \eqref{tsa0l} becomes an algebraic equation and we have the following incremental behaviour of the generators at time $t$ by linearizing \eqref{tsa0} with \eqref{tsa0l} reduced
\begin{equation}\label{tsaVaVb}
\begin{split}
   \begin{bmatrix}
     \bm{M}_a\Delta\ddot{\bm{\theta}}_a^t \\
     \bm{M}_b\Delta\ddot{\bm{\theta}}_b^t \\
   \end{bmatrix}+
   [\Lg(\bm{\theta}^t)/ \bm{L}_{cc}(\bm{\theta}^t)]
   \begin{bmatrix}
     \Delta\bm{\theta}_a^t \\
     \Delta\bm{\theta}_b^t \\
   \end{bmatrix}=\bm{0}_{a+b}
\end{split}
\end{equation}
where $\Delta\bm{\theta}^t=\bm{\theta}^{t+\Delta t}-\bm{\theta}^t$, $\Delta t\rightarrow0$ denotes the micro-increment of rotor angles after time $t$.
Note that \eqref{tsaVaVb} applies to any time of the post-fault system, which is different from \eqref{sda} that represents the linearization around an equilibrium point.

\indent
Since the generators in the same group have similar responses, their collective behaviours can be represented by the following two variables in the manner of center of inertia, say $\Delta\delta_a^t = m_a^{-1}\bm{1}_a^T\bm{M}_a\Delta\bm{\theta}_a^t$ and $\Delta\delta_b^t = m_b^{-1}\bm{1}_b^T\bm{M}_b\Delta\bm{\theta}_b^t$, where $m_a=\sum_{i\in \mV_a} M_i$ and $m_b=\sum_{i\in \mV_b} M_i$ denote the total inertia of $\mV_a$ and $\mV_b$, respectively.
By Petrov-Galerkin projection \cite{antoulas2005approximation}, we apply the approximation $\Delta\bm{\theta}_a^t=\Delta\delta_a^t\bm{1}_a$ and $\Delta\bm{\theta}_b^t=\Delta\delta_b^t\bm{1}_b$ so that \eqref{tsaVaVb} is reduced to
\begin{equation}\label{tsaVaVb2}
\begin{split}
   \Delta\ddot{\delta}_{ab}^t = -(m_a^{-1}+m_b^{-1})\geff(\Lg(\bm{\theta}^t), \mV_a, \mV_b)\Delta\delta_{ab}^t
\end{split}
\end{equation}
where $\Delta\delta_{ab}^t=\Delta\delta_a^t-\Delta\delta_b^t$.

\indent
If $\geff(\Lg(\bm{\theta}^t), \mV_a, \mV_b)$ becomes negative due to the disturbance, then the right hand side of \eqref{tsaVaVb2} functions as a positive feedback that tends to separate the rotor angles between the two groups of coherent generators. This tendency will further deteriorate the coupling between the two groups of generators and finally lead to transient instability as a group of generators losses synchronism with respect to the other one.
On the other hand, if we have a larger positive $\geff(\Lg(\bm{\theta}^e), \mV_a, \mV_b)$ at the operating equilibrium point, then the two groups of generators get more tightly coupled. Meanwhile, $\geff(\Lg(\bm{\theta}^t), \mV_a, \mV_b)$ can maintain positive and tend to prevent angle separation under more severe disturbances, which indicates higher stability.

\indent
Thus, given a fault and the corresponding groups of coherent generators $\mV_a, \mV_b$, we formulate a new transient stability constrained optimal power flow (TSC-OPF) model that makes the operating equilibrium point more robust to the fault
\begin{subequations}\label{OPFgeff}
\begin{align}
    \min_{P_i, \theta_i}~&\sum\nolimits_{i\in \mV_G} f_i(P_i)   \label{OPFgeffobj}\\
    s.t.~&P_i = \sum\nolimits_{j\in\mN_i} V_iV_jB_{ij}\sin\theta_{ij}, ~\forall i\in \mV \label{OPFgeffpf}\\
         & P_i^{\min} \leq P_i \leq P_i^{\max},~\forall i\in \mV \label{OPFgeffPmin}\\
         & V_i = V_i^{set},~\forall i\in \mV   \label{OPFgeffVmin}\\
         & |\theta_{ij}|\leq \theta_{ij}^{\max},~\forall (i,j)\in \mE \label{OPFgefftheta} \\
         &\geff(\Lg(\bm{\theta}), \mV_a, \mV_b)\geq g_{\min}   \label{OPFgeffgmin}
\end{align}
\end{subequations}
where the objective \eqref{OPFgeffobj} denotes the total generation cost that is usually a convex function of $P_i$;
\eqref{OPFgeffpf} refers to active power flow equation; \eqref{OPFgeffPmin} refers to power injection limit;
\eqref{OPFgeffVmin} regulates bus voltage magnitude;
\eqref{OPFgefftheta} refers to angle difference limit with $\theta_{ij}^{\max}<\frac{\pi}{2}$;
and $\geff(\Lg, \mV_a, \mV_b)$ is lower bounded by \eqref{OPFgeffgmin} to regulate post-fault system performance.
Problem \eqref{OPFgeff} provides new thoughts for transient stability enhancement, i.e., to strengthen the electrical coupling between the two groups of coherent generators via generation dispatch.

\indent
Problem \eqref{OPFgeff} is generally nonconvex due to the nonlinearity of \eqref{OPFgeffpf}, which adds difficulty to the solution method.
Nevertheless, convex relaxation technique \cite{low2014convex1, low2014convex2} provides a tractable way for OPF convexification, and problem \eqref{OPFgeff} admits a convex relaxation form due to the property of effective conductance.
To this end, we introduce the Hermitian matrix $\bm{U}=[U_{ij}]\in \mbC^{n\times n}$ where $U_{ii}=V_i^2$, $\forall i\in\mV$ and $U_{ij}=\overline{U}_{ji}=V_i V_j e^{\upj \theta_{ij}}$, $\forall (i,j)\in\mE$, and $U_{ij}=0$ otherwise.
In addition, we define the function $S_i(\bm{U}) = \sum\nolimits_{j\in\mN_i} \upj B_{ij}(U_{ii}-U_{ij})$.
Then, \eqref{OPFgeff} can be reformulated in terms of $\bm{U}$ as follows
\begin{subequations}\label{OPFgeffconvex}
\begin{align}
    \min_{P_i,\bm{U}}~&\sum\nolimits_{i\in \mV_G} f_i(P_i) \label{OPFgeffconvexobj}\\
    s.t.~&P_i = \textup{Re}\{S_i(\bm{U})\},~\forall i\in \mV    \label{OPFgeffconvexPequ}\\
          &P_i^{\min} \leq  P_i \leq P_i^{\max},~\forall i\in \mV    \label{OPFgeffconvexPmin}\\
         &U_{ii} = (V_i^{set})^2,~\forall i\in \mV   \label{OPFgeffconvexVmin}\\
      & |\textup{Im}\{U_{ij}\}|\leq\textup{Re}\{U_{ij}\}\tan\theta_{ij}^{\max},~\forall (i,j)\in \mE  \label{OPFgeffconvextheta1}\\
      &\begin{bmatrix}
         \bm{1}_a^T \bm{E}_a\bm{W}_q\bm{E}_a^T\bm{1}_a - g_{\min} & \bm{1}_a^T\bm{E}_a\bm{W}_q\bm{E}_c^T \\
         \bm{E}_c\bm{W}_q\bm{E}_a^T\bm{1}_a & \bm{E}_c\bm{W}_q\bm{E}_c^T \\
      \end{bmatrix}\succeq 0   \label{OPFgeffconvexgmin}\\
      &\bm{U}\succeq 0 \label{OPFgeffconvexSDP}\\
      &\textup{rank}(\bm{U}) = 1  \label{OPFgeffconvexrank}
\end{align}
\end{subequations}
where $\bm{W}_q(\bm{U})=\diag{B_{ij}\textup{Re}\{U_{ij}\}}\in \mbR^{l\times l}$, $\forall (i,j)\in \mE$.
The Laplacian matrix $\Lg=\bm{E}\bm{W}_q\bm{E}^T$ is PSD with only one zero eigenvalue as $\textup{Re}\{U_{ij}\}$ is restricted to be positive by \eqref{OPFgeffconvextheta1}. So constraint \eqref{OPFgeffconvexgmin} is equivalent to \eqref{OPFgeffgmin} following the same idea as in \eqref{geffschur}.
Then, according to \cite{low2014convex1}, problem \eqref{OPFgeffconvex} has the same feasible set as \eqref{OPFgeff}, and relaxing the rank constraint \eqref{OPFgeffconvexrank} gives a convex problem.
The convex relaxation of OPF has been well studied, which is shown to be exact in most cases \cite{low2014convex2}.
Hence, we can find the global optimum of the original problem \eqref{OPFgeff} by applying sophisticated solvers to the convex problem \eqref{OPFgeffconvexobj}-\eqref{OPFgeffconvexSDP}.
Even if the rank relaxation is inexact (i.e., $\textup{rank}(\bm{U})>1$), we can recover a rank-one solution from the convex problem by adding a very small penalty to the objective function that has little influence on the optimality~\cite{madani2015convex}.

\indent
Note that the proposed OPF model is mainly to confirm the potential application of effective conductance in transient stability enhancement.
The formulation is still preliminary and leaves much space for improvement.
For instance, the model only contains active power flow equation and stability constraint for a single fault.
A more comprehensive and practical model should include reactive power flow equation and stability constraints for a set of faults, which will be considered in future work.

\subsection{An illustrative case study}\label{seccase}
\begin{figure}[!h]
  \centering
  \includegraphics[width=1.5in]{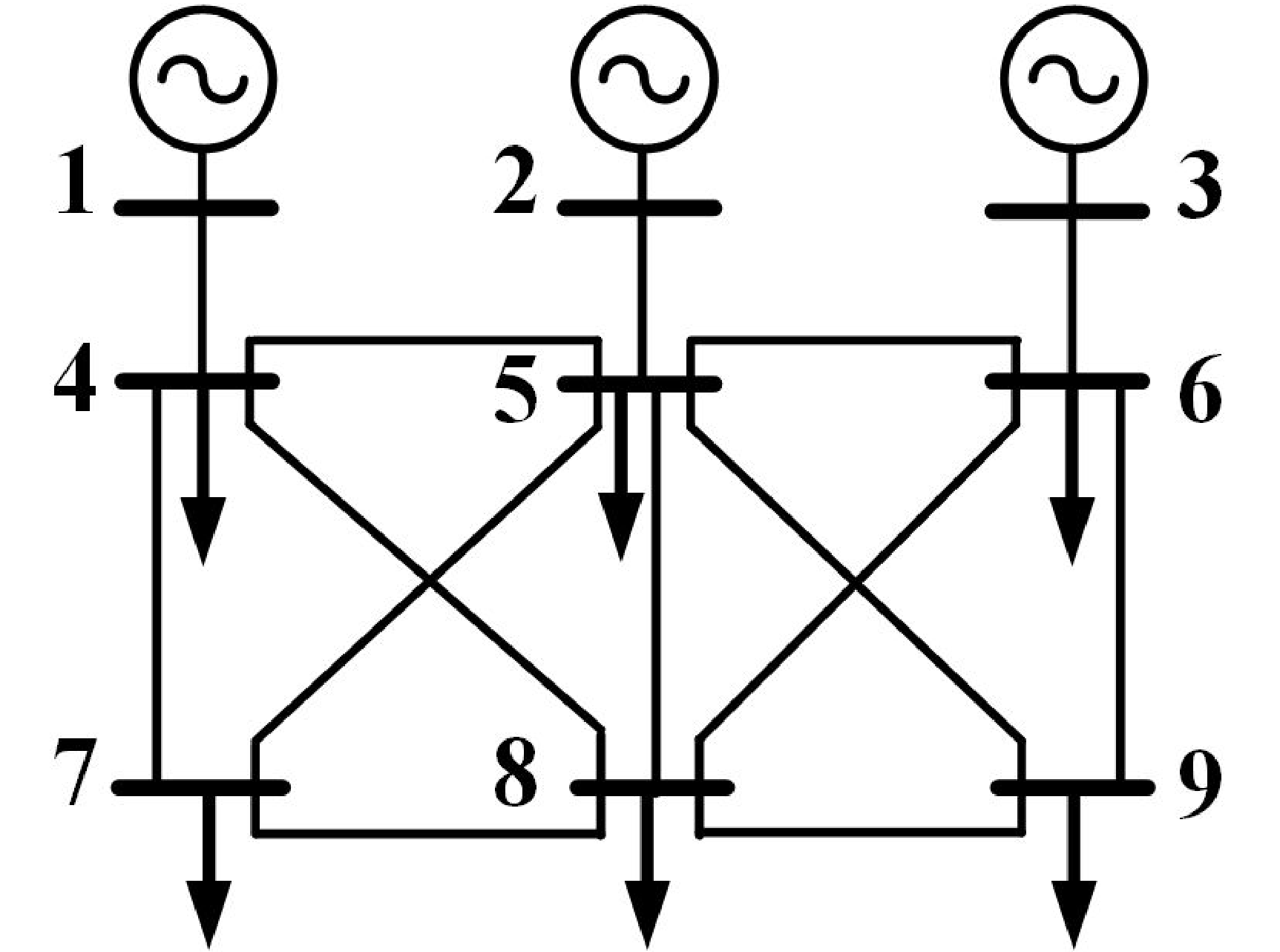}
  \caption{IEEE 6-bus system diagram (including generator internal buses).}
  \label{figcase6busdiagram}
\end{figure}

Take the IEEE 6-bus system to illustrate the obtained results.
The system diagram is shown in Fig. \ref{figcase6busdiagram}, where the network is augmented with three generator internal buses (bus 1, 2, 3) that respectively connect terminal bus 4, 5, 6. The system parameters (in per-unit value) are given in Table \ref{tab6buspara}.

\begin{table}[!h]
\renewcommand{\arraystretch}{1.3}
  \caption{IEEE 6-bus system parameters}
  \label{tab6buspara}
  \centering
    \begin{tabular}{cccc|cccc}
    \hline\hline
     Bus   &  $V_i^{set}$   &  $M_i$   &   $D_i$  &  Line  &  $B_{ij}$ &  Line  &  $B_{ij}$\\
    \hline
    1     & 1.05  &  10  &  2   & (4,5)     & 5.00   & (7,8)     & 2.50    \\
    2     & 1.05  &  20  &  2   & (4,7)     & 5.00   & (8,9)     & 3.33    \\
    3     & 1.07  &  30  &  2   & (4,8)     & 3.33   & (1,3)     & 50.00    \\
    4     & 1.05  &  --  &  2   & (5,6)     & 4.00    & (2,4)     & 50.00    \\
    5     & 1.05  &  --  &  2   & (5,7)     & 10.00   & (3,6)     & 50.00     \\
    6     & 1.07  &  --  &  2   & (5,8)     & 3.33     &  &    \\
    7     & 1.00  &  --  &  2   & (5,9)     & 5.00     &  &    \\
    8     & 1.00  &  --  &  2   & (6,8)     & 3.85    &  &    \\
    9     & 1.00  &  --  &  2   & (6,9)     & 10.00   &  &    \\
    \hline\hline
    \end{tabular}%
\end{table}%

\indent
We choose the equilibrium points under three power injection schemes in Table \ref{tab6bus3EPs} for small-disturbance stability analysis.
Equilibrium point A refers to a normal operating point where no critical lines exist, and it is trivial to check that $\Lg(\bm{\theta}^A)$ is PSD with only one zero eigenvalue and all eigenvalues of $\bm{J}_{dyn}(\bm{\theta}^A)$ have negative real parts.
There is one critical line (4, 8) at equilibrium point B, whose angle difference is 90.34\degree.
In this case, the effective conductance between bus 4 and bus 8 is 2.85, which implies that equilibrium point B is asymptotically stable by the obtained theorems.
On the other hand, direct calculation gives that $\Lg(\bm{\theta}^B)$ is PSD with only one zero eigenvalue and all eigenvalues of $\bm{J}_{dyn}(\bm{\theta}^B)$ have negative real parts, which coincides with the judgement given by the effective conductance.

\begin{figure}[!h]
  \centering
  \includegraphics[width=2.7in]{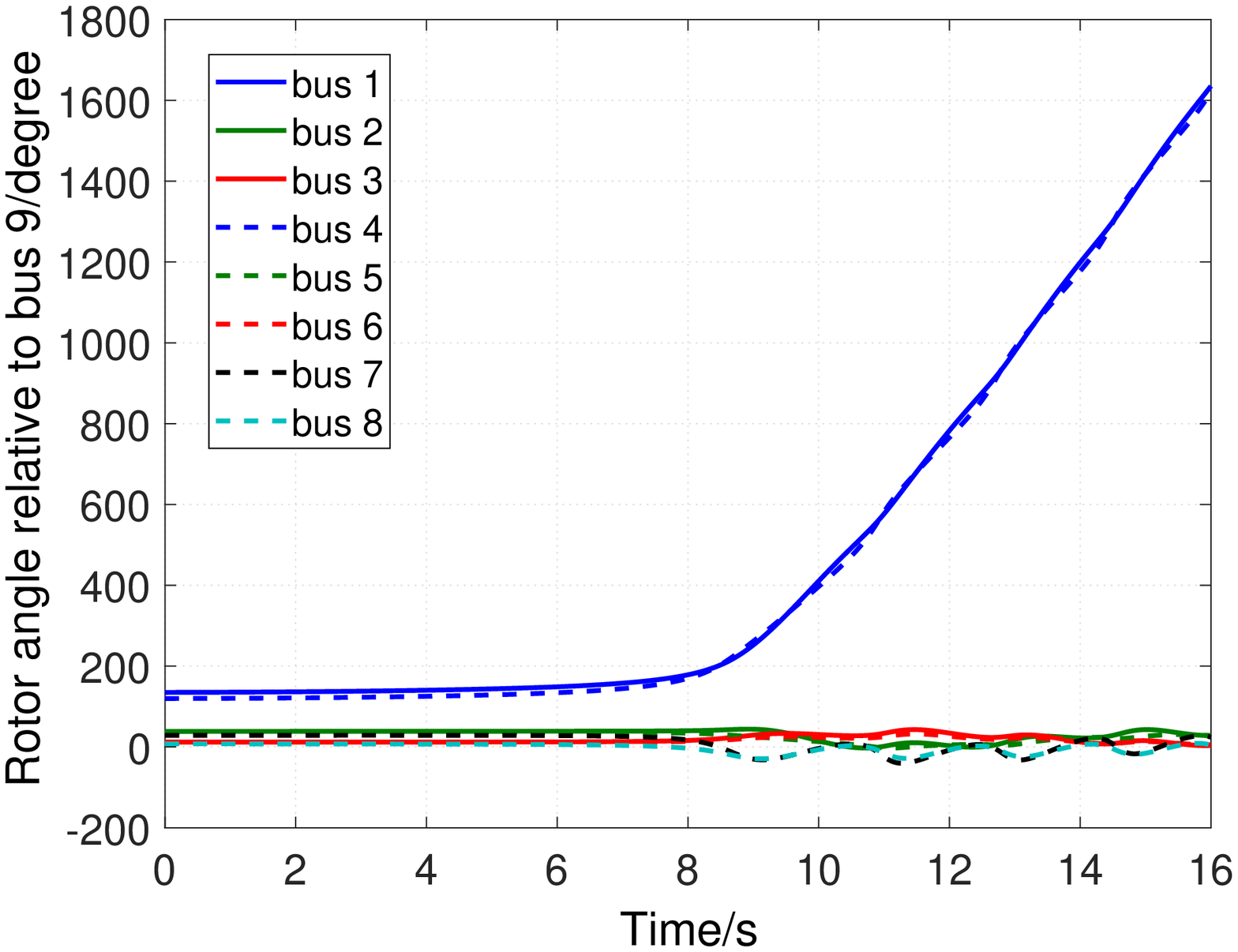}
  \caption{Angle curves after applying small disturbance to equilibrium point C.}
  \label{figsds}
\end{figure}

\indent
At equilibrium point C, there are two critical lines, namely line (4,7) and line (4,8), whose angle differences are 90.37\degree and 111.67\degree, respectively.
Negative effective conductances appear in this case, e.g., the effective conductance between bus 4 and bus 8 is -0.76.
By the obtained theorems, equilibrium point C is a UEP, $\Lg(\bm{\theta}^C)$ has at most two negative eigenvalues and $\bm{J}_{dyn}(\bm{\theta}^C)$ has at most two eigenvalues with positive real parts.
It coincides with the fact that $\Lg(\bm{\theta}^C)$ has a negative eigenvalue -0.550 and $\bm{J}_{dyn}(\bm{\theta}^C)$ has a positive eigenvalue 0.234.
Also, we find that the effective conductance hints the pattern of angle separation.
The effective conductance between bus 1 and bus set \{2,3\} is -0.91, and the effective conductance between bus 2 and bus 3 is 7.24.
So bus 2 and bus 3 are supposed to remain coherent and bus 1 tends to separate from bus 2, 3 when we apply a small disturbance to equilibrium point C. This speculation gets confirmed by the time domain simulation in Fig. \ref{figsds}.

\begin{figure}[!h]
  \centering
  \includegraphics[width=2.7in]{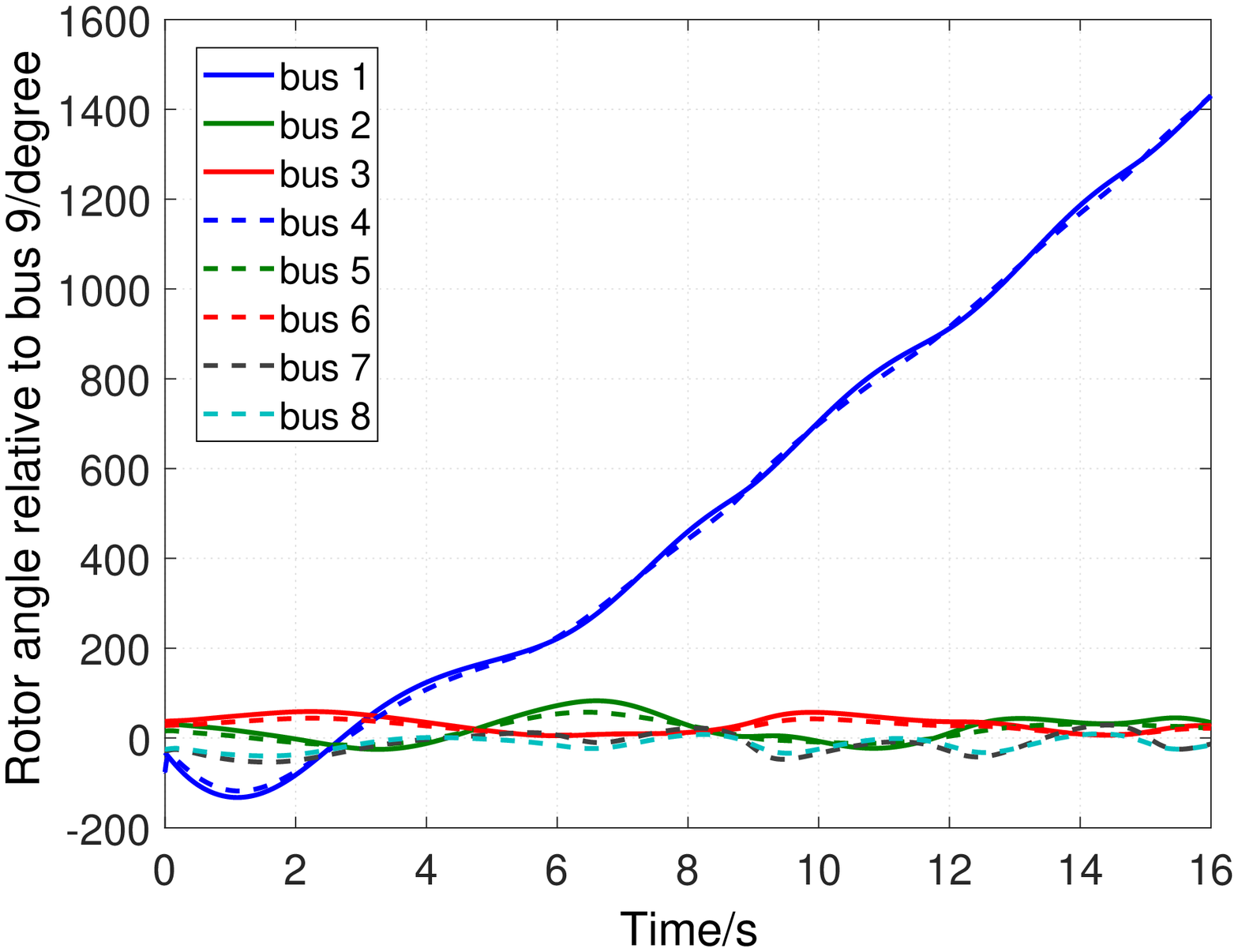}
  \caption{Angle curves after applying short-circuit fault to equilibrium point A with clearing time 1.0s.}
  \label{figfaultA}
\end{figure}

\begin{figure}[!h]
  \centering
  \includegraphics[width=2.7in]{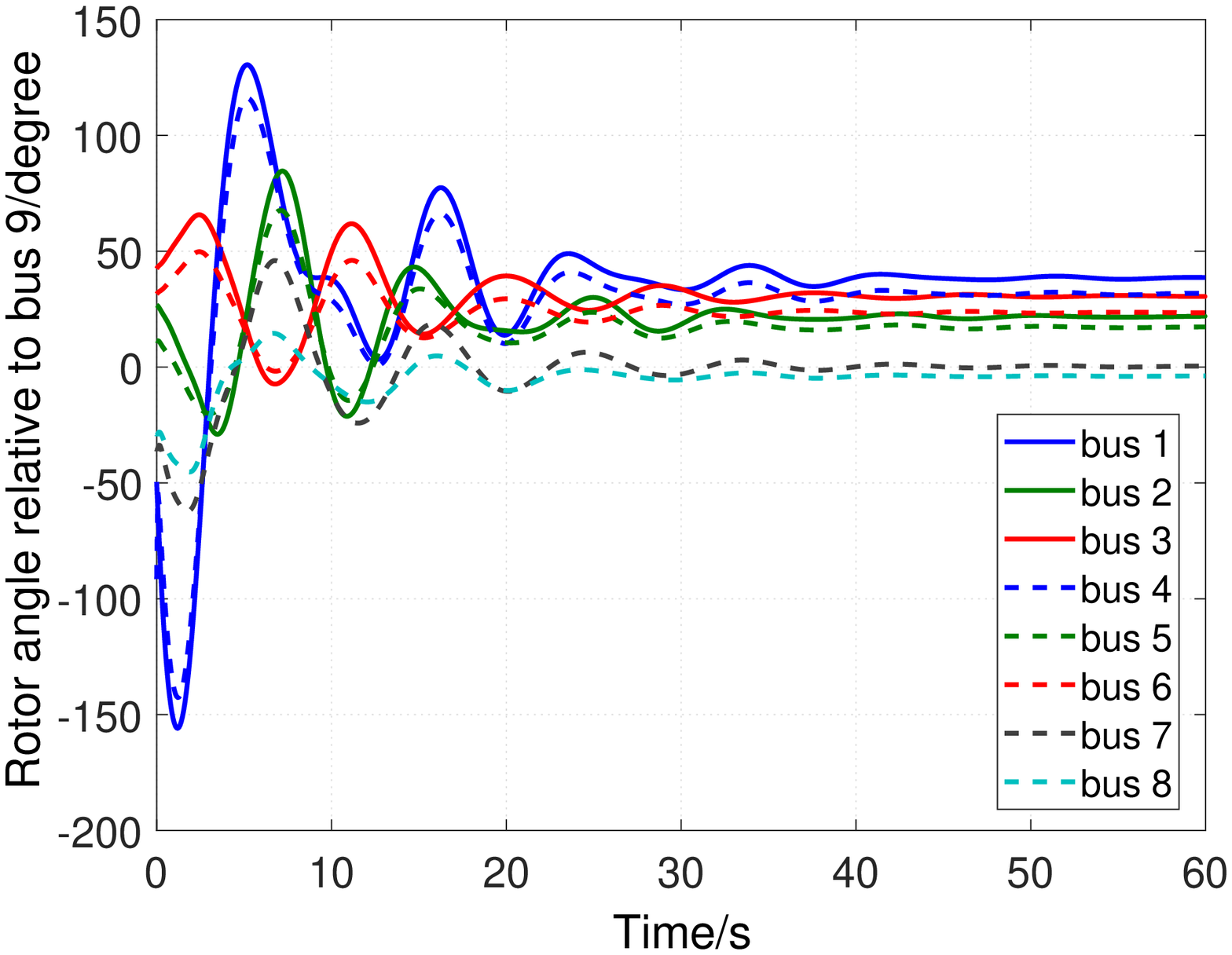}
  \caption{Angle curves after applying short-circuit fault to equilibrium point D with clearing time 1.0s.}
  \label{figfaultD}
\end{figure}

\indent
Next, we turn to the proposed TSC-OPF.
Suppose the system operates at equilibrium point A, and bus 4 (terminal bus of generator 1) has a three-phase short-circuit fault such that it is grounded via a 0.02 p.u. resistance.
If the fault is cleared after 1.0 second, the system is unstable as bus 1 and bus 4 lose synchronism with respect to the other buses (see Fig. \ref{figfaultA}).
In this case, the coherent groups of generators are $\mV_a=\{1\}$ and $\mV_b=\{2,3\}$, and the corresponding effective conductance at equilibrium point A is 7.58.
Then, we solve the relaxed convex problem of \eqref{OPFgeffconvex} with the settings $f_1=0.01 P_1^2+ 10P_1$, $f_2=0.01 P_2^2+ 12 P_2$, $f_3=0.01 P_3^2+ 12 P_3$; $P_i^{\min}=0$, $P_i^{\max}=15$, $\forall i\in\mV_G$; $P_i^{\min}=P_i^{\max}=P_i$, $\forall i\in\mV_L$; $\theta_{ij}^{\max}=45$\degree, $\forall (i,j)\in \mE$ and $g_{\min}=8$.
The program is implemented by using CVX toolbox and SeDuMi solver \cite{cvx}, and the obtained generation scheme (equilibrium point D) is shown in Table \ref{tab6busoptimization}. This scheme is globally optimal as the rank constraint \eqref{OPFgeffconvexrank} is numerically satisfied (the ratio of the largest eigenvalue to second largest eigenvalue of $\bm{U}$ is beyond $10^6$).
At equilibrium point D, the objective is 204.04 and the effective conductance between bus 1 and bus set \{2,3\} is increased to 8.00.
Consequently, the system achieves transient stability with the same fault and clearing time (see Fig. \ref{figfaultD}).

\indent
Another generation scheme (equilibrium point E) is listed in Table \ref{tab6busoptimization} for comparison, which is obtained by solving the OPF model without the effective conductance constraint \eqref{OPFgeffconvexgmin}.
At equilibrium point E, the objective is decreased to 198.88 and the effective conductance between bus 1 and bus set \{2,3\} is decreased to 7.30.
In this case, the system is unstable with the same fault and clearing time. The corresponding post-fault response is similar to Fig. \ref{figfaultA}, which is omitted here for simplicity.
This comparison shows the significance of including the effective conductance constraint into OPF model for the sake of transient stability.

\indent
Moreover, in Table \ref{tab6busCCT} we quantify the degree of transient stability of the five equilibrium points under this fault by using the critical clearing time (CCT), i.e., the maximum fault-on time that the system can suffer without losing stability.
We observe that the CCT increases with effective conductance, which further verifies that transient stability level grows with effective conductance.

\begin{table}[!h]
\renewcommand{\arraystretch}{1.3}
  \caption{Three equilibrium points \& corresponding generation schemes}
  \label{tab6bus3EPs}
  \centering
    \begin{tabular}{c|cc|cc|cc}
    \hline\hline
    Bus   & $P_i^A$    & $\theta_i^A$ & $P_i^B$    & $\theta_i^B$ & $P_i^C$    & $\theta_i^C$  \\
    \hline
    1     & 8.40  &  0.00\degree   & 13.80 & 0.00\degree  & 14.50 & 0.00\degree   \\
    2     & 4.40  & -27.86\degree & 3.40  & -74.33\degree & 2.70  & -96.16\degree  \\
    3     & 5.20  & -29.86\degree & 0.80  & -101.40\degree & 0.80  & -122.97\degree  \\
    4     & -0.40 & -8.76\degree   & -0.40 & -14.50\degree & -0.40 & -15.25\degree  \\
    5     & -0.40 & -32.44\degree & -0.40 & -77.87\degree & -0.40 & -98.96\degree  \\
    6     & -0.40 & -35.07\degree & -0.40 & -102.20\degree & -0.40 & -123.77\degree  \\
    7     & -5.60 & -47.19\degree & -5.60 & -85.78\degree & -5.60 & -105.62\degree  \\
    8     & -5.60 & -54.91\degree & -5.60 & -104.83\degree & -5.60 & -126.92\degree  \\
    9     & -5.60 & -54.79\degree & -5.60 & -113.32\degree & -5.60 & -134.87\degree  \\
    \hline\hline
    \end{tabular}%
\end{table}%

\begin{table}[!h]
\renewcommand{\arraystretch}{1.3}
  \caption{The operating point after optimization}
  \label{tab6busoptimization}
  \centering
    \begin{tabular}{c|cc|cc}
    \hline\hline
    Bus   & $P_i^D$   & $\theta_i^D$ &  $P_i^E$   & $\theta_i^E$ \\
    \hline
    1     &   6.54  &    0.00\degree &  9.37  &    0.00\degree  \\
    2     &   4.39  & -16.69\degree &  0.00  &   -41.55\degree \\
    3     &   7.07  & -7.79\degree &    8.63  &   -16.00\degree \\
    4     &   -0.40 & -6.81\degree &  -0.40  &   -9.78\degree \\
    5     &  -0.40 & -21.25\degree &  -0.40 &    -41.55\degree  \\
    6     &  -0.40 & -14.89\degree &  -0.40  &   -24.67\degree \\
    7     &  -5.60 & -38.15\degree &  -5.60  &   -53.00\degree \\
    8     &  -5.60 & -42.34\degree &  -5.60  &   -54.79\degree \\
    9     &  -5.60 & -38.44\degree &  -5.60  &   -51.68\degree \\
    \hline\hline
    \end{tabular}%
\end{table}%

\begin{table}[!h]
\renewcommand{\arraystretch}{1.3}
  \caption{Effective conductance \& CCT for fault at bus 4}
  \label{tab6busCCT}
  \centering
    \begin{tabular}{ccc}
    \hline\hline
    Generation scheme   & $\geff(\Lg, \{1\}, \{2,3\})$    & CCT (s) \\
    \hline
    Equilibrium point C     &   -0.91  & 0.00 (UEP)  \\
    Equilibrium point B     &   3.47  & 0.22  \\
    Equilibrium point E     &   7.30  & 0.87   \\
    Equilibrium point A     &   7.58  & 0.95   \\
    Equilibrium point D     &   8.00  & 1.12  \\
    \hline\hline
    \end{tabular}%
\end{table}%

\section{Conclusion and future direction}\label{secconclu}
The effective resistance and effective conductance have been extended from node-to-node versions to nodes-to-nodes ones.
The proposed definitions properly characterize the relation between any two disjoint sets of nodes in the network with clear circuit interpretations. The convexity and monotonicity of the proposed definitions have been studied.
Then, these new definitions lead to an intuitive interpretation of graph Laplacian definiteness.
It has been proved that graph Laplacian being PSD with only one zero eigenvalue is equivalent to the effective conductances between specific pairs of node sets being positive. The number of Laplacian negative eigenvalues is upper bounded by the number of negative weighted edges.
The effective conductance is also a powerful indicator for power network stability.
Necessary and sufficient conditions in terms of effective conductance has been established for the small-disturbance angle stability, hyperbolicity and type of power system equilibria.
In addition, a TSC-OPF model has been developed based on the effective conductance, which provides a new way for transient stability enhancement. The proposed TSC-OPF can be reformulated into a convex relaxation form due to properties of effective conductance, and hence achieves global optimality.

\indent
Future works include developing a more comprehensive TSC-OPF model with more practical factors included, such as reactive power flow and effective conductance constraints for a fault set.
Also the effective resistance/conductance may be of interest to other issues in power networks such as energy management~\cite{liu2019high}, which remains to be explored.

\section*{Appendix}
\subsection{Proofs in Section \ref{secdefreff}}
   \textit{Proof of Lemma \ref{lemLcc}: }
   It is trivial that
   \begin{subequations}
   \begin{align}
     \bm{0}_a &= \bm{L}_{aa}\bm{1}_{a} + \bm{L}_{ab}\bm{1}_{b} + \bm{L}_{ac}\bm{1}_{c} \label{LG11}\\
     \bm{0}_b &= \bm{L}_{ab}^T\bm{1}_{a} + \bm{L}_{bb}\bm{1}_{b} + \bm{L}_{bc}\bm{1}_{c} \label{LG12}\\
     \bm{0}_c &= \bm{L}_{ac}^T\bm{1}_{a} + \bm{L}_{bc}^T\bm{1}_{b} + \bm{L}_{cc}\bm{1}_{c} \label{LG13}
   \end{align}
   \end{subequations}
   We have $\bm{1}_{c}=-\bm{L}_{cc}^{-1}(\bm{L}_{ac}^T\bm{1}_{a} + \bm{L}_{bc}^T\bm{1}_{b})$ from \eqref{LG13}
   Then, substituting it into \eqref{LG11} and \eqref{LG12} leads to \eqref{Lcc1}.

   \indent
   By Definition \ref{defreff} and $(\Lg/ \bm{L}_{cc})\bm{1}_{a+b}=\bm{0}_{a+b}$, we have
   \begin{equation*}
   \begin{split}
   &\geff(\Lg, \mV_b, \mV_a) =
   \begin{bmatrix}
        \bm{0}_a^T  &   \bm{1}_{b}^T  \\
     \end{bmatrix}
     (\Lg/ \bm{L}_{cc})
     \begin{bmatrix}
     \bm{0}_a \\
     \bm{1}_b  \\
   \end{bmatrix}\\
   &=-\begin{bmatrix}
        \bm{1}_a^T  &   \bm{0}_b^T  \\
     \end{bmatrix}
     (\Lg/ \bm{L}_{cc})
     \begin{bmatrix}
     \bm{0}_a \\
     \bm{1}_b  \\
   \end{bmatrix}\\
   &= \begin{bmatrix}
        \bm{1}_a^T  &   \bm{0}_b^T  \\
     \end{bmatrix}
     (\Lg/ \bm{L}_{cc})
     \begin{bmatrix}
     \bm{1}_a \\
     \bm{0}_b  \\
   \end{bmatrix}=\geff(\Lg, \mV_a, \mV_b)
   \end{split}
   \end{equation*}
   which implies that \eqref{Lcc2} holds.

   \indent
   We expand $\Lg/ \bm{L}_{cc}$ as follows
    \begin{equation*}
    \begin{split}
     \Lg/ \bm{L}_{cc} =
     \begin{bmatrix}
     \bm{L}_{aa} - \bm{L}_{ac}\bm{L}_{cc}^{-1}\bm{L}_{ac}^T   &   \bm{L}_{ab} - \bm{L}_{ac}\bm{L}_{cc}^{-1}\bm{L}_{bc}^T    \\
     \bm{L}_{ab}^T - \bm{L}_{bc}\bm{L}_{cc}^{-1}\bm{L}_{ac}^T   &   \bm{L}_{bb} - \bm{L}_{bc}\bm{L}_{cc}^{-1}\bm{L}_{bc}^T   \\
     \end{bmatrix}.
    \end{split}
    \end{equation*}
    Then \eqref{Lcc3aa}-\eqref{Lcc3bb} follow from $(\Lg/ \bm{L}_{cc})\bm{1}_{a+b}=\bm{0}_{a+b}$.   \hfill $\blacksquare$

    \textit{Proof of Lemma \ref{lemquotient}: }
    For $\mV_r\subset\mV$, $\mV_o=\mV\backslash\mV_r$ and two disjoint non-empty sets $\mV_a, \mV_b\subset \mV_r$, let $\mV_p=\mV_r\backslash (\mV_a\cup\mV_b)$, $\mV_q=\mV_a\cup\mV_b$ so that $\mV_c=\mV\backslash(\mV_a\cup\mV_b)=\mV_p\cup\mV_o$.
    Then, $\mV$ is divided into three disjoint sets, say $\mV_o$, $\mV_p$ and $\mV_q$, and we can expand the Laplacian matrix of $\mG_r(\mV_r, \mE_r, \bm{W}_r)$ as
    \begin{equation*}
    \begin{split}
     \Lg/ \bm{L}_{oo} =
     \begin{bmatrix}
     \bm{L}_{pp} - \bm{L}_{po}\bm{L}_{oo}^{-1}\bm{L}_{po}^T   &   \bm{L}_{pq} - \bm{L}_{po}\bm{L}_{oo}^{-1}\bm{L}_{qo}^T    \\
     \bm{L}_{pq}^T - \bm{L}_{qo}\bm{L}_{oo}^{-1}\bm{L}_{po}^T   &   \bm{L}_{qq} - \bm{L}_{qo}\bm{L}_{oo}^{-1}\bm{L}_{qo}^T   \\
     \end{bmatrix}
    \end{split}
    \end{equation*}
    where the block $\bm{L}_{pp} - \bm{L}_{po}\bm{L}_{oo}^{-1}\bm{L}_{po}^T$ refers to the submatrix whose rows and columns are indexed by $\mV_p$.
    Further, we observe that $\bm{L}_{pp} - \bm{L}_{po}\bm{L}_{oo}^{-1}\bm{L}_{po}^T=\bm{L}_{cc}/ \bm{L}_{oo}$ where $\bm{L}_{cc} =
     \begin{bmatrix}
     \bm{L}_{pp}   &   \bm{L}_{po}    \\
     \bm{L}_{po}^T   &   \bm{L}_{oo}  \\
     \end{bmatrix}$.
    Since $\bm{L}_{cc}$ is assumed nonsingular, it follows from Lemma \ref{lemschur} that $\bm{L}_{cc}/\bm{L}_{oo}$ is nonsingular.
    Then, by Definition \ref{defreff} we have
    \begin{equation*}
    \begin{split}
      \geff(\Lgr, \mV_a, \mV_b)=\bm{e}_{\mV_a}^T[(\Lg/ \bm{L}_{oo})/(\bm{L}_{cc}/ \bm{L}_{oo})]\bm{e}_{\mV_a}.
    \end{split}
    \end{equation*}
    The quotient formula of Schur complement \cite{zhang2006schur} gives $(\Lg/ \bm{L}_{oo})/(\bm{L}_{cc}/ \bm{L}_{oo})=\Lg/ \bm{L}_{cc}$, so we conclude that $\geff(\Lgr, \mV_a, \mV_b)=\bm{e}_{\mV_a}^T(\Lg/ \bm{L}_{cc})\bm{e}_{\mV_a}=\geff(\Lg, \mV_a, \mV_b)$.
    \hfill $\blacksquare$

   \textit{Proof of Theorem \ref{thmenergyVdiff}: }
   1) First, we prove that \eqref{energydis} holds.
   We transform the optimization problem in \eqref{energydis} into the following equivalent form by the relation $\bm{i}=\Lg\bm{v}$
   \begin{subequations}\label{energydisv}
   \begin{align}
      \min_{\bm{v}}~&\bm{v}^T\Lg\bm{v} \label{energydis0v} \\
      s.t.~&\bm{e}_{\mV_a}^T\Lg\bm{v} = 1  \label{energydis1v}\\
        &\bm{e}_{\mV_b}^T\Lg\bm{v} = -1  \label{energydis2v}\\
        &\bm{e}_i^T\Lg\bm{v} = 0,~\forall i\in \mV_c. \label{energydis3v}
   \end{align}
   \end{subequations}
   Note that the objective function \eqref{energydis0v} is derived from
   \begin{equation*}
   \begin{split}
        \bm{i}^T\Lg^{\dag}\bm{i} = \bm{v}^T\Lg\Lg^{\dag}\Lg\bm{v}
        =\bm{v}^T\Lg(\bm{I}_n-\frac{1}{n}\bm{1}_n\bm{1}_n^T)\bm{v} = \bm{v}^T\Lg\bm{v}
   \end{split}
   \end{equation*}
   by applying Lemma \ref{lemLdag}.
   Also we point out that problem \eqref{energydisv} is convex since $\Lg\succeq 0$.
   Thus, the global optimum of problem \eqref{energydisv} is given by the following KKT condition \cite{boyd2004convex}
   \begin{subequations}\label{KKT}
   \begin{align}
        &\bm{e}_{\mV_a}^T\Lg\bm{v} = 1 \\
        &\bm{e}_{\mV_b}^T\Lg\bm{v} = -1  \\
        &\bm{e}_i^T\Lg\bm{v} = 0,~i\in \mV_c\\
        &2\Lg\bm{v}+
        \begin{bmatrix}
          \bm{L}_{aa}\bm{1}_a \\
          \bm{L}_{ab}^T\bm{1}_a \\
          \bm{L}_{ac}^T\bm{1}_a \\
        \end{bmatrix}\mu_a+
        \begin{bmatrix}
          \bm{L}_{ab}\bm{1}_b \\
          \bm{L}_{bb}\bm{1}_b \\
          \bm{L}_{bc}^T\bm{1}_b \\
        \end{bmatrix}\mu_b+
        \begin{bmatrix}
          \bm{L}_{ac} \\
          \bm{L}_{bc} \\
          \bm{L}_{cc} \\
        \end{bmatrix}\bm{\mu}_c  \label{KKT4}
        =\bm{0}_n
   \end{align}
   \end{subequations}
   where $\mu_a\in\mbR, \mu_b\in\mbR, \bm{\mu}_c\in\mbR^c$ are the Lagrangian multipliers corresponding to constraint \eqref{energydis1v}, \eqref{energydis2v}, \eqref{energydis3v}, respectively.
   Let $\mu_a^{*} = -2\reff(\Lg, \mV_a, \mV_b)$, $\mu_b^{*} = 0$, $\bm{\mu}_c^{*} = 2\reff(\Lg, \mV_a, \mV_b)\bm{L}_{cc}^{-1}\bm{L}_{ac}^T\bm{1}_a$ and
   \begin{equation}
   \begin{split}
        \bm{v}^{*} = \reff(\Lg, \mV_a, \mV_b)\Lg^{\dag}
        \begin{bmatrix}
         (\bm{L}_{aa}-\bm{L}_{ac}\bm{L}_{cc}^{-1}\bm{L}_{ac}^T)\bm{1}_a\\
         (\bm{L}_{ab}^T-\bm{L}_{bc}\bm{L}_{cc}^{-1}\bm{L}_{ac}^T)\bm{1}_a\\
         \bm{0}_c\\
        \end{bmatrix}.
   \end{split}
   \end{equation}
   Since \eqref{Lcc3aa} and \eqref{Lcc3ab} lead to $\bm{1}_a^T(\bm{L}_{aa}-\bm{L}_{ac}\bm{L}_{cc}^{-1}\bm{L}_{ac}^T)\bm{1}_a+
   \bm{1}_b^T(\bm{L}_{ab}^T-\bm{L}_{bc}\bm{L}_{cc}^{-1}\bm{L}_{ac}^T)\bm{1}_a=0$, we have
   \begin{equation}\label{Lvprime}
   \begin{split}
       &\Lg\bm{v}^{*} =\\
       & \reff(\Lg, \mV_a, \mV_b)(\bm{I}_n-\frac{1}{n}\bm{1}_n\bm{1}_n^T)
       \begin{bmatrix}
         (\bm{L}_{aa}-\bm{L}_{ac}\bm{L}_{cc}^{-1}\bm{L}_{ac}^T)\bm{1}_a\\
         (\bm{L}_{ab}^T-\bm{L}_{bc}\bm{L}_{cc}^{-1}\bm{L}_{ac}^T)\bm{1}_a\\
         \bm{0}_c\\
        \end{bmatrix}\\&=\reff(\Lg, \mV_a, \mV_b)
        \begin{bmatrix}
         (\bm{L}_{aa}-\bm{L}_{ac}\bm{L}_{cc}^{-1}\bm{L}_{ac}^T)\bm{1}_a\\
         (\bm{L}_{ab}^T-\bm{L}_{bc}\bm{L}_{cc}^{-1}\bm{L}_{ac}^T)\bm{1}_a\\
         \bm{0}_c\\
        \end{bmatrix}.
   \end{split}
   \end{equation}
   Thus, it can be seen that $(\bm{v}^{*}, \mu_a^{*}, \mu_b^{*}, \bm{\mu}_c^{*})$ satisfies the KKT condition \eqref{KKT} and hence $\bm{v}^{*}$ is the global optimum of \eqref{energydisv}.

   \indent
   Observing that \eqref{KKT4} can be rewritten as
   $\Lg(2\bm{v}^{*}+
        \begin{bmatrix}
          \bm{1}_a^T\mu_a^{*} &  \bm{0}_b^T &  (\bm{\mu}_c^{*})^T \\
        \end{bmatrix}^T)=\bm{0}_n$,
   which implies that
   \begin{equation}\label{vprime}
   \begin{split}
       (\bm{v}^{*})^T = k\bm{1}_n^T-\frac{1}{2}
       \begin{bmatrix}
          (\bm{1}_a\mu_a^{*})^T & \bm{0}_b^T & (\bm{\mu}_c^{*})^T \\
        \end{bmatrix},~k\in \mbR
   \end{split}
   \end{equation}
   since $\Lg$ has only one zero eigenvalue with the eigenvector being $\bm{1}_n$.
   Then, substituting \eqref{Lvprime} and \eqref{vprime} into \eqref{energydis0v} gives
   \begin{equation}\label{vLv}
   \begin{split}
       &(\bm{v}^{*})^T\Lg\bm{v}^{*}= \\
       &k\cdot\reff(\Lg, \mV_a, \mV_b)\bm{1}_n^T
        \begin{bmatrix}
         (\bm{L}_{aa}-\bm{L}_{ac}\bm{L}_{cc}^{-1}\bm{L}_{ac}^T)\bm{1}_a\\
         (\bm{L}_{ab}^T-\bm{L}_{bc}\bm{L}_{cc}^{-1}\bm{L}_{ac}^T)\bm{1}_a\\
         \bm{0}_c\\
        \end{bmatrix}-\\
        &\frac{\reff(\Lg, \mV_a, \mV_b)}{2}
        \begin{bmatrix}
          \bm{1}_a\mu_a^{*} \\
         \bm{0}_b \\
          \bm{\mu}_c^{*} \\
        \end{bmatrix}^T
       \begin{bmatrix}
         (\bm{L}_{aa}-\bm{L}_{ac}\bm{L}_{cc}^{-1}\bm{L}_{ac}^T)\bm{1}_a\\
         (\bm{L}_{ab}^T-\bm{L}_{bc}\bm{L}_{cc}^{-1}\bm{L}_{ac}^T)\bm{1}_a\\
         \bm{0}_c\\
        \end{bmatrix}.
   \end{split}
   \end{equation}
   By \eqref{Lcc3aa} and \eqref{Lcc3ab}, the first term of right-hand-side of \eqref{vLv} is zero, and the second term is $\reff(\Lg, \mV_a, \mV_b)$.
   Therefore, $\reff(\Lg, \mV_a, \mV_b)$ is the optimal value of problem \eqref{energydisv} as well as the original problem in \eqref{energydis}.
   Moreover, we have $\reff(\Lg, \mV_a, \mV_b)>0$ since the objective $\bm{v}^T\Lg\bm{v}$ is non-negative and the only case that makes it zero ($\bm{v}=\bm{1}_n$) is excluded by \eqref{energydis1v}.

   \indent
   2) Second, we prove that \eqref{Vdiff} holds.
   Note that constraints \eqref{energydis1v}-\eqref{energydis3v} can be expressed by the compact form $\bm{\mP}_{ab}^T\Lg\bm{v}=\bm{e}_1-\bm{e}_2\in\mbR^{c+2}$.
   Let us regard \eqref{energydisv} as a primal problem, then its Lagrange dual problem takes the form
   \begin{equation*}
   \begin{split}
      f(\bm{v},\bm{\mu})=\bm{v}^T\Lg\bm{v} + \bm{\mu}^T(\bm{\mP}_{ab}^T\Lg\bm{v}-\bm{e}_1+\bm{e}_2)
   \end{split}
   \end{equation*}
   where $\bm{\mu}\in \mbR^{c+2}$ is the Lagrange multiplier vector.
   Since problem \eqref{energydisv} is convex and strictly feasible, we have the strong duality below \cite{boyd2004convex}
   \begin{equation}\label{energydisdual}
   \begin{split}
      \reff(\Lg, \mV_a, \mV_b)=\max_{\bm{\mu}}~\min_{\bm{v}}~f(\bm{v},\bm{\mu}).
   \end{split}
   \end{equation}

   \indent
   Next, we derive the explicit expression for the optimal value of problem \eqref{energydisdual}.
   The internal problem $\min_{\bm{v}}~f(\bm{v},\bm{\mu})$ is convex with respect to $\bm{v}$, and the (global) optimum is given by
   \begin{equation}\label{zeroderivative}
   \begin{split}
       \frac{\partial f}{\partial \bm{v}}=2\Lg\bm{v}+\Lg\bm{\mP}_{ab}\bm{\mu}=\bm{0}_n.
   \end{split}
   \end{equation}
   Since $\Lg$ has only one zero eigenvalue, \eqref{zeroderivative} implies that the optimum is $\bm{v}=-\frac{1}{2}\bm{\mP}_{ab}\bm{\mu}+k\bm{1}_n$, $k\in \mbR$, so that problem \eqref{energydisdual} is reduced to
   \begin{equation}\label{energydisdual2}
   \begin{split}
      \reff(\Lg, \mV_a, \mV_b)=\max_{\bm{\mu}}~-\frac{1}{4}\bm{\mu}^T\bm{\mP}_{ab}^T\Lg\bm{\mP}_{ab}\bm{\mu} -\bm{\mu}^T(\bm{e}_1-\bm{e}_2).
   \end{split}
   \end{equation}
   Problem \eqref{energydisdual2} is convex with respect to $\bm{\mu}$, and the (global) optimum is given by
   \begin{equation}\label{zeroderivative2}
   \begin{split}
       \frac{1}{2}\bm{\mP}_{ab}^T\Lg\bm{\mP}_{ab}\bm{\mu}+\bm{e}_1-\bm{e}_2=\bm{0}_{c+2}.
   \end{split}
   \end{equation}
   Note that $\bm{\mP}_{ab}^T\Lg\bm{\mP}_{ab}$ is the Laplacian matrix of the graph with $\mV_a$ and $\mV_b$ being clustered, so $\bm{\mP}_{ab}^T\Lg\bm{\mP}_{ab}$ has only one zero eigenvalue and it follows from Lemma \ref{lemLdag} that
   \begin{equation}\label{Ldag}
   \begin{split}
       (\bm{\mP}_{ab}^T\Lg\bm{\mP}_{ab})^{\dag}\bm{\mP}_{ab}^T\Lg\bm{\mP}_{ab} = \bm{I}_{c+2}-\frac{1}{c+2}\bm{1}_{c+2}\bm{1}_{c+2}^T.
   \end{split}
   \end{equation}
   Then, by \eqref{zeroderivative2} and \eqref{Ldag}, the optimum of \eqref{energydisdual2} is
   \begin{equation}\label{optimalmu}
   \begin{split}
       \bm{\mu}=-2(\bm{\mP}_{ab}^T\Lg\bm{\mP}_{ab})^{\dag}(\bm{e}_1-\bm{e}_2)+k\bm{1}_{c+2},~k\in \mbR.
   \end{split}
   \end{equation}
   Finally, substituting \eqref{optimalmu} into \eqref{energydisdual2} gives $\reff(\Lg, \mV_a, \mV_b)=(\bm{e}_1-\bm{e}_2)^T(\bm{\mP}_{ab}^T\Lg\bm{\mP}_{ab})^{\dag}(\bm{e}_1-\bm{e}_2)$. \hfill $\blacksquare$

\subsection{Proofs in Section \ref{secmonotone}}
   \textit{Proof of Theorem \ref{thmmonotone}: }
   We first prove that $\geff(\Lg, \mV_a, \mV_b)$ is non-decreasing by discussing the following situations for $\Delta w_{ij}$:

   \indent
   1) If $i,j\in \mV_a$, then $\Delta w_{ij}$ only influences $\bm{L}_{aa}$.
   Without loss of generality, we suppose $i,j$ are the first and second node in $\mV_a$, then by \eqref{Lcc3aa} we have
   \begin{equation*}
   \begin{split}
       \Delta\geff(\Lg, \mV_a, \mV_b) = \Delta w_{ij}\bm{1}_{a}^T(\bm{e}_1-\bm{e}_2)(\bm{e}_1-\bm{e}_2)^T\bm{1}_{a} = 0.
   \end{split}
   \end{equation*}
   Similarly, we have the same result if $i,j\in \mV_b$.

   \indent
   2) If $i,j\in \mV_c$, then $\Delta w_{ij}$ only influences $\bm{L}_{cc}$.
   Without loss of generality, we suppose $i,j$ are the first and second node in $\mV_c$.
   By \eqref{Lcc3aa} and Sherman-Morrison formula [\citen{petersen2008matrix}, p.18], we have
   \begin{equation*}
   \begin{split}
      &\Delta \geff(\Lg, \mV_a, \mV_b) =\\ &-\bm{1}_{a}^T\bm{L}_{ac}(\bm{L}_{cc}+\Delta\bm{L}_{cc})^{-1}\bm{L}_{ac}^T\bm{1}_{a} +\bm{1}_{a}^T\bm{L}_{ac}\bm{L}_{cc}^{-1}\bm{L}_{ac}^T\bm{1}_{a}\\
      &= \bm{1}_{a}^T\bm{L}_{ac}\frac{\Delta w_{ij}\bm{L}_{cc}^{-1}(\bm{e}_{1}-\bm{e}_{2})(\bm{e}_{1}-\bm{e}_{2})^T\bm{L}_{cc}^{-1}}
      {1+\Delta w_{ij}(\bm{e}_{1}-\bm{e}_{2})^T\bm{L}_{cc}^{-1}(\bm{e}_{1}-\bm{e}_{2})}\bm{L}_{ac}^T\bm{1}_{a}.
   \end{split}
   \end{equation*}
   Note that $\bm{L}_{cc}$ is positive definite as $\Lg$ is PSD with only one zero eigenvalue \cite{horn2012matrix}, and hence $\Delta \geff(\Lg, \mV_a, \mV_b)\geq 0$.

   \indent
   3) If $i\in \mV_a$, $j\in \mV_b$, then $\Delta w_{ij}$ only influences the matrices $\bm{L}_{aa}$, $\bm{L}_{ab}$, $\bm{L}_{bb}$.
   Without loss of generality, we suppose $i,j$ are respectively the first node in $\mV_a$ and $\mV_b$, by \eqref{Lcc3aa} we have
   \begin{equation*}
   \begin{split}
      \Delta\geff(\Lg, \mV_a, \mV_b) = \Delta w_{ij}\bm{1}_{a}^T\bm{e}_1\bm{e}_1^T\bm{1}_{a} > 0.
   \end{split}
   \end{equation*}
   Similarly, we have the same result if $i\in \mV_b$, $j\in \mV_a$.

   \indent
   4) If $i\in \mV_a$, $j\in \mV_c$, then $\Delta w_{ij}$ only influences the matrices $\bm{L}_{aa}, \bm{L}_{cc}, \bm{L}_{ac}$.
   Without loss of generality, we suppose $i,j$ are respectively the first node in $\mV_a$ and $\mV_c$.
   Since $\bm{L}_{cc}$ is positive definite, by \eqref{Lcc3bb} and Sherman-Morrison formula we have
   \begin{equation*}
   \begin{split}
      \Delta \geff(\Lg, \mV_a, \mV_b) = \bm{1}_{b}^T\bm{L}_{bc}\frac{\Delta w_{ij}\bm{L}_{cc}^{-1}\bm{e}_{1}\bm{e}_{1}^T\bm{L}_{cc}^{-1}}
      {1+\Delta w_{ij}\bm{e}_{1}^T\bm{L}_{cc}^{-1}\bm{e}_{1}}\bm{L}_{bc}^T\bm{1}_{b}\geq 0.
   \end{split}
   \end{equation*}
   Similarly, we have the same result if $i\in \mV_c$, $j\in \mV_a$ or $i\in \mV_c$, $j\in \mV_b$ or $i\in \mV_b$, $j\in \mV_c$.
   So we have proved that $\geff(\Lg, \mV_a, \mV_b)$ is non-decreasing after adding an edge $(i,j)$ with weight $\Delta w_{ij}>0$ in all possible cases.

   \indent
   In addition, since $\reff(\Lg, \mV_a, \mV_b)=1/\geff(\Lg, \mV_a, \mV_b)>0$ if $\Lg$ is PSD with only one zero eigenvalue, we conclude that $\reff(\Lg, \mV_a, \mV_b)$ is non-increasing after adding an edge $(i,j)$ with weight $\Delta w_{ij}>0$.   \hfill $\blacksquare$

\subsection{Proofs in Section \ref{seclapPSD}}
   \textit{Proof of Lemma \ref{lemPSD}: }
   Without loss of generality, the sets in the sequential inclusion of $\mV$ are numbered as $\mV^{n-k}=\{k+1,...,n\}$, $k=0,...,n-1$.
   Then condition \eqref{geffpos} becomes $\geff(\Lg, \{k\}, \mV^{n-k})>0$, $k=1,...,n-1$.
   Thus, for each $k=1,...,n-1$, we have the corresponding $\mV_a=\{k\}$, $\mV_b=\mV^{n-k}$, $\mV_c=\{1,...,k-1\}$, and $\bm{L}_{cc}$ is nonsingular (we set $\mV_c=\phi$ and ignore $\bm{L}_{cc}$ when $k=1$).
   It implies that $[\Lg]_k$, $k=1,...,n-2$ is nonsingular, where $[\Lg]_k$ denotes the $k$-th order leading principal submatrix of $\Lg$.

   \indent
   We now show by contradiction that $[\Lg]_{n-1}$ is also nonsingular.
   Suppose $[\Lg]_{n-1}$ is singular with $[\Lg]_{n-1}\bm{x}=\bm{0}_{n-1}$, where $\bm{x}\in\mbR^{n-1}$ is a nonzero vector.
   Express $\Lg$ as $\Lg=
      \begin{bmatrix}
         [\Lg]_{n-1} & \bm{a} \\
         \bm{a}^T & L_{nn} \\
       \end{bmatrix}$, where $\bm{a}\in\mbR^{n-1}$ is indexed by the first $n-1$ rows and last column of $\Lg$.
   Then we have
   \begin{equation*}
   \begin{split}
      \bm{x}^T\bm{a}&=\bm{x}^T\bm{a}+\bm{x}^T[\Lg]_{n-1}\bm{1}_{n-1}\\
      &=\bm{x}^T([\Lg]_{n-1}\bm{1}_{n-1}+\bm{a})=\bm{x}^T\bm{0}_{n-1}=0
   \end{split}
   \end{equation*}
   where $[\Lg]_{n-1}\bm{1}_{n-1}+\bm{a}=\bm{0}_{n-1}$ is given by the zero row sums of $\Lg$.
   It follows that
   $\Lg\begin{bmatrix}
         \bm{x} \\
         0 \\
      \end{bmatrix}=
      \begin{bmatrix}
           [\Lg]_{n-1}\bm{x} \\
           \bm{a}^T\bm{x} \\
      \end{bmatrix}=\bm{0}_n$, which implies that $\Lg$ has a second zero eigenvalue with the eigenvector being
   $\begin{bmatrix}
         \bm{x}^T &  0 \\
      \end{bmatrix}^T$.
   On the other hand, consider the case when $\mV_a=\{n-1\}$, $\mV_b=\mV^{1}$, $\mV_c=\{1,...,n-2\}$.
   In this case, we have $\geff(\Lg, \{n-1\}, \mV^{1})>0$ and $\Lg/\bm{L}_{cc}$ is a two-by-two matrix.
   Then, it follows from Lemma \ref{lemLcc} that
   $\Lg/\bm{L}_{cc} = \geff(\Lg, \{n-1\}, \mV^{1})\cdot
      \begin{bmatrix}
         1 & -1 \\
         -1 & 1 \\
       \end{bmatrix}$, which gives $i_0(\Lg/\bm{L}_{cc})=1$.
   Further, we have $i_0(\Lg)=i_0(\Lg/\bm{L}_{cc})=1$ by Lemma \ref{lemschur}.
   It implies that $\Lg$ has only one zero eigenvalue with the eigenvector being $\bm{1}_n$, which yields a contradiction.
   Therefore, we conclude that $[\Lg]_k$ is nonsingular for $k=1,...,n-1$.

   \indent
   Similar to \eqref{currentflow}, we have the following circuit equation by letting $\mV_a=\{k\}$, $\mV_b=\mV^{n-k}$ and $\mV_c=\{1,...,k-1\}$
   \begin{equation*}
   \begin{split}
      \begin{bmatrix}
      \bm{i}_c^T  &   i_a & \bm{i}_b^T  \\
     \end{bmatrix}^T=
      \begin{bmatrix}
      \bm{0}_c^T  &   i_k  &  \bm{i}_b^T
     \end{bmatrix}^T= \Lg
      \begin{bmatrix}
      \bm{v}_c^T   &   1   &  \bm{0}_b^T \\
     \end{bmatrix}^T
   \end{split}
   \end{equation*}
   the first two rows of which can be rewritten as
   \begin{equation}\label{currentVac}
   \begin{split}
      \begin{bmatrix}
      \bm{i}_c   \\
       i_a  \\
     \end{bmatrix}=
      \begin{bmatrix}
      \bm{0}_c   \\
      i_k  \\
     \end{bmatrix}= [\Lg]_k
      \begin{bmatrix}
      \bm{v}_c   \\
      1   \\
     \end{bmatrix}.
   \end{split}
   \end{equation}
   Since $[\Lg]_k$, $k=1,...,n-1$ is nonsingular, by Theorem \ref{thmIa} and \eqref{currentVac} we have
   \begin{equation}\label{geffI10}
   \begin{split}
      &\geff(\Lg, \{k\}, \mV^{n-k}) =  i_k =
      \begin{bmatrix}
      \bm{0}_c^T   &   i_k  \\
     \end{bmatrix}
     \begin{bmatrix}
      \bm{v}_c   \\
      1   \\
     \end{bmatrix}\\ &=
     \begin{bmatrix}
      \bm{0}_c^T    &   i_k \\
     \end{bmatrix} [\Lg]_k^{-1}
      \begin{bmatrix}
      \bm{0}_c   \\
      i_k   \\
     \end{bmatrix},~k=1,...,n-1
   \end{split}
   \end{equation}
   which implies that $\geff(\Lg, \{k\}, \mV^{n-k})$ equals to the product of $i_k^2$ and $(k,k)$-entry of $[\Lg]_k^{-1}$.
   Expressing the matrix inverse $[\Lg]_k^{-1}$ in \eqref{geffI10} by its cofactors [\citen{petersen2008matrix}, p.17] gives
   \begin{equation}\label{geffI1}
   \begin{split}
      \geff(\Lg, \{k\}, \mV^{n-k}) =
      \left\{
        \begin{array}{ll}
           i_k^2\frac{1}{\textup{det}([\Lg]_{1})},~k=1 \\
           i_k^2\frac{\textup{det}([\Lg]_{k-1})}{\textup{det}([\Lg]_k)},~k=2,...,n-1.
        \end{array}
      \right.
   \end{split}
   \end{equation}
   Applying $\geff(\Lg, \{1\}, \mV^{n-1})>0$ to \eqref{geffI1} gives $\textup{det}([\Lg]_{1})>0$.
   Similarly, iteratively applying $\geff(\Lg, \{k\}, \mV^{n-k})>0$, $k=2,...,n-1$ to \eqref{geffI1} gives $\textup{det}([\Lg]_k)>0$, $k=2,...,n-1$.
   Thus, the matrix $[\Lg]_{n-1}$ is positive definite as all its leading principal minors of order $k$, $k\leq n-1$ are positive \cite{horn2012matrix}.
   Then, by Lemma \ref{lemschur} we have $i_+(\Lg)\geq i_+([\Lg]_{n-1})=n-1$, and hence $\Lg$ is PSD with only one zero eigenvalue.   \hfill $\blacksquare$

\ifCLASSOPTIONcaptionsoff
  \newpage
\fi

{\footnotesize
\bibliographystyle{IEEEtran}
\bibliography{effresistance}

\begin{thebibliography}{10}
\providecommand{\url}[1]{#1}
\csname url@samestyle\endcsname
\providecommand{\newblock}{\relax}
\providecommand{\bibinfo}[2]{#2}
\providecommand{\BIBentrySTDinterwordspacing}{\spaceskip=0pt\relax}
\providecommand{\BIBentryALTinterwordstretchfactor}{4}
\providecommand{\BIBentryALTinterwordspacing}{\spaceskip=\fontdimen2\font plus
\BIBentryALTinterwordstretchfactor\fontdimen3\font minus
  \fontdimen4\font\relax}
\providecommand{\BIBforeignlanguage}[2]{{%
\expandafter\ifx\csname l@#1\endcsname\relax
\typeout{** WARNING: IEEEtran.bst: No hyphenation pattern has been}%
\typeout{** loaded for the language `#1'. Using the pattern for}%
\typeout{** the default language instead.}%
\else
\language=\csname l@#1\endcsname
\fi
#2}}
\providecommand{\BIBdecl}{\relax}
\BIBdecl

\bibitem{klein1993resistance}
D.~J. Klein and M.~Randi{\'c}, ``Resistance distance,'' \emph{J. Math. Chem.},
  vol.~12, no.~1, pp. 81--95, 1993.

\bibitem{foster1949average}
R.~M. Foster, ``The average impedance of an electrical network,'' in
  \emph{Contributions to Applied Mechanics (Reissner Anniversary
  Volume)}.\hskip 1em plus 0.5em minus 0.4em\relax Ann Arbor, MI, USA: Edwards
  Brothers, 1949, pp. 333--340.

\bibitem{bapat2010graphs}
R.~B. Bapat, \emph{Graphs and Matrices}.\hskip 1em plus 0.5em minus 0.4em\relax
  Springer, 2010.

\bibitem{ghosh2008minimizing}
A.~Ghosh, S.~Boyd, and A.~Saberi, ``Minimizing effective resistance of a
  graph,'' \emph{SIAM Review}, vol.~50, no.~1, pp. 37--66, 2008.

\bibitem{thulasiraman2019network}
K.~Thulasiraman, M.~Yadav, and K.~Naik, ``Network science meets circuit theory:
  Resistance distance, {K}irchhoff index, and {F}oster's theorems with
  generalizations and unification,'' \emph{{IEEE} Trans. Circuits Syst. {I}},
  vol.~66, no.~3, pp. 1090--1103, Mar. 2019.

\bibitem{coppersmith1993random}
D.~Coppersmith, P.~Doyle, P.~Raghavan \emph{et~al.}, ``Random walks on weighted
  graphs and applications to on-line algorithms,'' \emph{Journal of the ACM},
  vol.~40, no.~3, pp. 421--453, 1993.

\bibitem{bonchev1994molecular}
D.~Bonchev, A.~T. Balaban, X.~Liu \emph{et~al.}, ``Molecular cyclicity and
  centricity of polycyclic graphs. i. cyclicity based on resistance distances
  or reciprocal distances,'' \emph{International Journal of Quantum Chemistry},
  vol.~50, no.~1, pp. 1--20, 1994.

\bibitem{guo2017monotonicity}
L.~Guo, C.~Liang, and S.~H. Low, ``Monotonicity properties and spectral
  characterization of power redistribution in cascading failures,'' in
  \emph{Annual Allerton Conference}, 2017, pp. 918--925.

\bibitem{soltan2017analysis}
S.~Soltan, D.~Mazauric, and G.~Zussman, ``Analysis of failures in power
  grids,'' \emph{{IEEE} Trans. Control Netw. Syst.}, vol.~4, no.~2, pp.
  288--300, Jun. 2017.

\bibitem{cotilla2013multi}
E.~Cotilla-Sanchez, P.~D. Hines, C.~Barrows \emph{et~al.}, ``Multi-attribute
  partitioning of power networks based on electrical distance,'' \emph{{IEEE}
  Trans. Power Syst.}, vol.~28, no.~4, pp. 4979--4987, Nov. 2013.

\bibitem{dorfler2010spectral}
F.~Dorfler and F.~Bullo, ``Spectral analysis of synchronization in a lossless
  structure-preserving power network model,'' in \emph{IEEE SmartGridComm.},
  2010, pp. 179--184.

\bibitem{dorfler2018electrical}
F.~D{\"o}rfler, J.~W. Simpson-Porco, and F.~Bullo, ``Electrical networks and
  algebraic graph theory: Models, properties, and applications,'' \emph{Proc.
  {IEEE}}, vol. 106, no.~5, pp. 977--1005, May 2018.

\bibitem{song2017networkbased}
Y.~Song, D.~J. Hill, and T.~Liu, ``Network-based analysis of small-disturbance
  angle stability of power systems,'' \emph{{IEEE} Trans. Control Netw. Syst.},
  vol.~5, no.~3, pp. 901--912, Sep. 2018.

\bibitem{barooah2006graph}
P.~Barooah and J.~P. Hespanha, ``Graph effective resistance and distributed
  control: Spectral properties and applications,'' in \emph{Proc. {IEEE} Conf.
  Dec. Control}, 2006, pp. 3479--3485.

\bibitem{grunberg2018performance}
T.~W. Grunberg and D.~F. Gayme, ``Performance measures for linear oscillator
  networks over arbitrary graphs,'' \emph{{IEEE} Trans. Control Netw. Syst.},
  vol.~5, no.~1, pp. 456--468, Mar. 2018.

\bibitem{johansson2018optimization}
A.~Johansson, J.~Wei, H.~Sandberg \emph{et~al.}, ``Optimization of the
  $\mathcal{H}_ {\infty}$-norm of dynamic flow networks,'' in \emph{Proc. Amer.
  Control Conf.}, 2018, pp. 1280--1285.

\bibitem{young2016new1}
G.~F. Young, L.~Scardovi, and N.~E. Leonard, ``A new notion of effective
  resistance for directed graphs---{P}art {I}: Definition and properties,''
  \emph{{IEEE} Trans. Autom. Control}, vol.~61, no.~7, pp. 1727--1736, Jul.
  2016.

\bibitem{young2016new2}
G.~F. Young, L.~Scardovi, and N.~E. Leonard, ``A new notion of effective
  resistance for directed graphs---{P}art {II}: Computing resistances,''
  \emph{{IEEE} Trans. Autom. Control}, vol.~61, no.~7, pp. 1737--1752, Jul.
  2016.

\bibitem{chen2016characterizing}
W.~Chen, J.~Liu, Y.~Chen \emph{et~al.}, ``Characterizing the positive
  semidefiniteness of signed {L}aplacians via effective resistances,'' in
  \emph{Proc. {IEEE} Conf. Dec. Control}, 2016, pp. 985--990.

\bibitem{zelazo2014definiteness}
D.~Zelazo and M.~B{\"u}rger, ``On the definiteness of the weighted {L}aplacian
  and its connection to effective resistance,'' in \emph{Proc. {IEEE} Conf.
  Dec. Control}, Dec 2014, pp. 2895--2900.

\bibitem{zelazo2017robustness}
D.~Zelazo and M.~B{\"u}rger, ``On the robustness of uncertain consensus
  networks,'' \emph{{IEEE} Trans. Control Netw. Syst.}, vol.~4, no.~2, pp.
  170--178, 2017.

\bibitem{chen2016definiteness}
Y.~Chen, S.~Z. Khong, and T.~T. Georgiou, ``On the definiteness of graph
  {L}aplacians with negative weights: Geometrical and passivity-based
  approaches,'' in \emph{Proc. Amer. Control Conf.}, 2016, pp. 2488--2493.

\bibitem{ahmadizadeh2017eigenvalues}
S.~Ahmadizadeh, I.~Shames, S.~Martin \emph{et~al.}, ``On eigenvalues of
  {L}aplacian matrix for a class of directed signed graphs,'' \emph{Linear
  Algebra and its Applications}, vol. 523, pp. 281--306, 2017.

\bibitem{chen2017spectral}
W.~Chen, D.~Wang, J.~Liu \emph{et~al.}, ``On spectral properties of signed
  {L}aplacians for undirected graphs,'' in \emph{Proc. {IEEE} Conf. Dec.
  Control}, 2017, pp. 1999--2002.

\bibitem{song2017cutset}
Y.~Song, D.~J. Hill, and T.~Liu, ``Characterization of cutsets in networks with
  application to transient stability analysis of power systems,'' \emph{{IEEE}
  Trans. Control Netw. Syst.}, vol.~5, no.~3, pp. 1261--1274, Sep. 2018.

\bibitem{gutman2004generalized}
I.~Gutman and W.~Xiao, ``Generalized inverse of the {L}aplacian matrix and some
  applications.'' \emph{Bulletin. Classe des Sciences Math{\'e}matiques et
  Naturelles. Sciences Math{\'e}matiques}, vol. 129, no.~29, pp. 15--23, 2004.

\bibitem{zhang2006schur}
F.~Zhang, \emph{The Schur Complement and Its Applications}.\hskip 1em plus
  0.5em minus 0.4em\relax Springer Science \& Business Media, 2006.

\bibitem{doyle1984random}
P.~G. Doyle and J.~L. Snell, \emph{Random Walks and Electric Networks}.\hskip
  1em plus 0.5em minus 0.4em\relax Washington, DC: Math. Assoc. America, 1984.

\bibitem{dorfler2013kron}
F.~D{\"o}rfler and F.~Bullo, ``Kron reduction of graphs with applications to
  electrical networks,'' \emph{{IEEE} Trans. Circuits Syst. {I}}, vol.~60,
  no.~1, pp. 150--163, Jan. 2013.

\bibitem{chua1980dynamic}
L.~O. Chua, ``Dynamic nonlinear networks: State-of-the-art,'' \emph{{IEEE}
  Trans. Circuits Syst.}, vol.~27, no.~11, pp. 1059--1087, Nov. 1980.

\bibitem{van2010characterization}
A.~Van~der Schaft, ``Characterization and partial synthesis of the behavior of
  resistive circuits at their terminals,'' \emph{Systems \& Control Letters},
  vol.~59, no.~7, pp. 423--428, 2010.

\bibitem{chua1987linear}
L.~O. Chua, C.~A. Desoer, and E.~S. Kuh, \emph{Linear and Nonlinear
  Circuits}.\hskip 1em plus 0.5em minus 0.4em\relax McGraw-Hill College, 1987.

\bibitem{altafini2013consensus}
C.~Altafini, ``Consensus problems on networks with antagonistic interactions,''
  \emph{{IEEE} Trans. Autom. Control}, vol.~58, no.~4, pp. 935--946, April
  2013.

\bibitem{boyd2004convex}
S.~Boyd and L.~Vandenberghe, \emph{Convex Optimization}.\hskip 1em plus 0.5em
  minus 0.4em\relax Cambridge University Press, 2004.

\bibitem{pan2016laplacian}
L.~Pan, H.~Shao, and M.~Mesbahi, ``Laplacian dynamics on signed networks,'' in
  \emph{Proc. {IEEE} Conf. Dec. Control}, 2016, pp. 891--896.

\bibitem{kundur2004definition}
P.~Kundur, J.~Paserba, V.~Ajjarapu \emph{et~al.}, ``Definition and
  classification of power system stability,'' \emph{{IEEE} Trans. Power Syst.},
  vol.~19, no.~3, pp. 1387--1401, Aug 2004.

\bibitem{bergen1981structure}
A.~R. Bergen and D.~J. Hill, ``A structure preserving model for power system
  stability analysis,'' \emph{{IEEE} Trans. Power App. Syst.}, vol. PAS-100,
  no.~1, pp. 25--35, Jan 1981.

\bibitem{chiang1987foundations}
H.-D. Chiang, F.~F. Wu, and P.~Varaiya, ``Foundations of direct methods for
  power system transient stability analysis,'' \emph{{IEEE} Trans. Circuits
  Syst.}, vol.~34, no.~2, pp. 160--173, Feb. 1987.

\bibitem{horn2012matrix}
R.~A. Horn and C.~R. Johnson, \emph{Matrix Analysis}.\hskip 1em plus 0.5em
  minus 0.4em\relax Cambridge University Press, 2012.

\bibitem{wu1983identification}
F.~F. Wu and Y.-K. Tsai, ``Identification of groups of $\epsilon$-coherent
  generators,'' \emph{{IEEE} Trans. Circuits Syst.}, vol.~30, no.~4, pp.
  234--241, Apr. 1983.

\bibitem{vu2016toward}
T.~L. Vu, S.~M. Al~Araifi, M.~S. El~Moursi \emph{et~al.}, ``Toward
  simulation-free estimation of critical clearing time,'' \emph{{IEEE} Trans.
  Power Syst.}, vol.~31, no.~6, pp. 4722--4731, Nov. 2016.

\bibitem{kundur1994power}
P.~Kundur, \emph{Power System Stability and Control}.\hskip 1em plus 0.5em
  minus 0.4em\relax McGraw-hill New York, 1994.

\bibitem{antoulas2005approximation}
A.~C. Antoulas, \emph{Approximation of Large-Scale Dynamical Systems}.\hskip
  1em plus 0.5em minus 0.4em\relax SIAM, 2005.

\bibitem{low2014convex1}
S.~H. Low, ``Convex relaxation of optimal power flow---{P}art {I}: Formulations
  and equivalence,'' \emph{{IEEE} Trans. Control Netw. Syst.}, vol.~1, no.~1,
  pp. 15--27, Mar. 2014.

\bibitem{low2014convex2}
S.~H. Low, ``Convex relaxation of optimal power flow---{P}art {II}:
  Exactness,'' \emph{{IEEE} Trans. Control Netw. Syst.}, vol.~1, no.~2, pp.
  177--189, Jun. 2014.

\bibitem{madani2015convex}
R.~Madani, S.~Sojoudi, and J.~Lavaei, ``Convex relaxation for optimal power
  flow problem: Mesh networks,'' \emph{{IEEE} Trans. Power Syst.}, vol.~30,
  no.~1, pp. 199--211, Jan. 2015.

\bibitem{cvx}
M.~Grant and S.~Boyd, ``{CVX}: Matlab software for disciplined convex
  programming, version 2.1,'' \url{http://cvxr.com/cvx}, Mar. 2014.

\bibitem{liu2019high}
G.~Liu, C.~Yuan, X.~Chen \emph{et~al.}, ``A high-performance energy management
  system based on evolving graph,'' \emph{{IEEE} Trans. Circuits Syst. {II}},
  2019, available online.

\bibitem{petersen2008matrix}
K.~B. Petersen and M.~S. Pedersen, \emph{The Matrix Cookbook}.\hskip 1em plus
  0.5em minus 0.4em\relax Technical University of Denmark, 2008.

\end{thebibliography}

}




\end{document}